# COMPUTATIONAL LIMITS OF A DISTRIBUTED ALGORITHM FOR SMOOTHING SPLINE

By Zuofeng Shang[†] and Guang Cheng[*]

*Indiana University-Purdue University at Indianapolis and Purdue University*

July 19, 2017

In this paper, we explore statistical versus computational trade-off to address a basic question in the application of a distributed algorithm: what is the minimal computational cost in obtaining statistical optimality? In smoothing spline setup, we observe a phase transition phenomenon for the number of deployed machines that ends up being a simple proxy for computing cost. Specifically, a *sharp* upper bound for the number of machines is established: when the number is below this bound, statistical optimality (in terms of nonparametric estimation or testing) is achievable; otherwise, statistical optimality becomes impossible. These sharp bounds partly capture intrinsic computational limits of the distributed algorithm considered in this paper, and turn out to be fully determined by the smoothness of the regression function. As a side remark, we argue that sample splitting may be viewed as an alternative form of regularization, playing a similar role as smoothing parameter.

**1. Introduction.** In the parallel computing environment, divide-and-conquer (D&C) method distributes data to multiple machines, and then aggregates local estimates computed from each machine to produce a global one. Such a distributed algorithm often requires a *growing* number of machines in order to process an increasingly large dataset. A practically relevant question is "how many processors do we really need in this parallel computing?" or "shall we allocate all our computational resources in the data analysis?" Such questions are related to the minimal computational cost of this distributed method (which will be defined more precisely later).

The major goal of this paper is to provide some "theoretical" insights for the above questions from a statistical perspective. Specifically, we consider a classical nonparametric regression setup:

$$(1.1) \qquad y_l = f(l/N) + \epsilon_l, \ l = 0, 1, \ldots, N-1,$$

where $\epsilon_l$'s are *iid* random errors with $E\{\epsilon_l\} = 0$ and $Var(\epsilon_l) = 1$, in the following distributed

---

[*]Assistant Professor.

[†]Corresponding Author. Professor. Research Sponsored by NSF (CAREER Award DMS-1151692, DMS-1418042), Simons Fellowship in Mathematics and Office of Naval Research (ONR N00014-15-1-2331). Guang Cheng gratefully acknowledges Statistical and Applied Mathematical Sciences Institute (SAMSI) for the hospitality and support during his visit in the 2013-Massive Data Program.

algorithm:

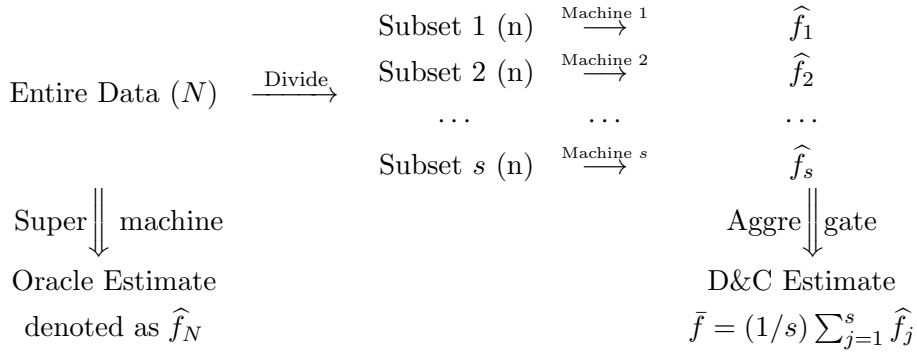

We assume that the total sample size is $N$, the number of machines is $s$ and the size of each sub-sample is $n$. Hence, $N = s \times n$. Each machine produces an individual smoothing spline estimate $\widehat{f}_j$ to be defined in (2.2) ([13]).

A known property of the above D&C strategy is that it can preserve statistical efficiency for a wide-ranging choice of $s$ (as demonstrated in Figure 1), say $\log s / \log N \in [0, 0.4]$, while largely reducing computational burden as $\log s / \log N$ increases (as demonstrated in Figure 2). An important observation from Figure 1 is that there is an obvious blowup for mean squared errors of $\bar{f}$ when the above ratio is beyond some threshold, e.g, 0.8 for $N = 10000$. Hence, we are interested in knowing whether there exists a critical value of $\log s / \log N$ in theory, beyond which statistical optimality no longer exists. For example, mean squared errors will never achieve minimax optimal lower bound (at rate level) no matter how smoothing parameters are tuned. Such a sharpness result partly captures the computational limit of the particular D&C algorithm considered in this paper, also complementing the upper bound results in [10, 16, 17]

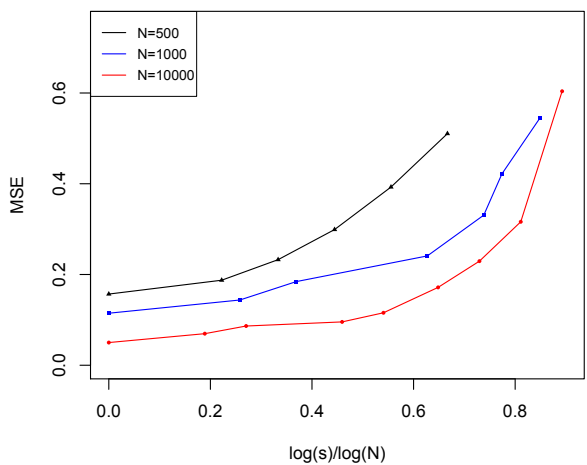

FIG 1. *Mean-square errors (MSE) of $\bar{f}$ based on 500 independent replications under different choices of $N$ and $s$. The values of MSE stay at low levels for various choice of $s$ with $\log s / \log N \in [0, 0.7]$. True regression function is $f_0(z) = 0.6b_{30,17}(z) + 0.4b_{3,11}(z)$ with $b_{a_1,a_2}$ the density function for $Beta(a_1, a_2)$.*

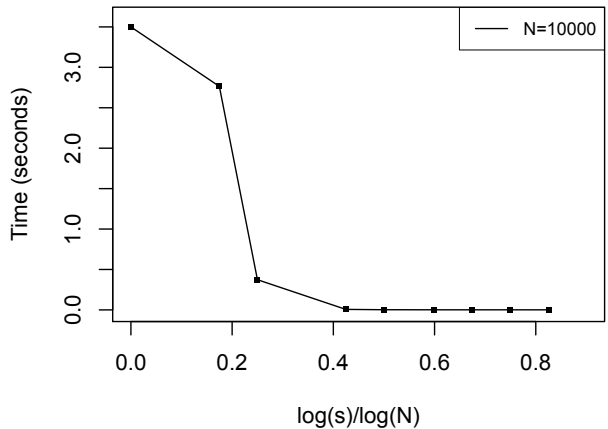

FIG 2. *Computing time of $\bar{f}$ based on a single replication under different choices of $s$ when $N = 10,000$. The larger the $s$, the smaller the computing time.*

Our first contribution is to establish a sharp upper bound of $s$ under which $\bar{f}$ achieves the minimax optimal rate $N^{m/(2m+1)}$, where $m$ represents the smoothness of $f_0$. By "sharp" upper bound, we mean the largest possible upper bound for $s$ to gain statistical optimality. This result is established by directly computing (non-asymptotic) upper and lower bounds of mean squared error of $\bar{f}$. These two bounds hold *uniformly* as $s$ diverges, and thus imply that the rate of mean squared error transits once $s$ reaches the rate $N^{2m/(2m+1)}$, which we call as phase transition in divide-and-conquer estimation. In fact, the choice of smoothing parameter, denoted as $\lambda$, also plays a very subtle role in the above phase transition. For example, $\lambda$ is not necessarily chosen at an optimal level when $s$ attains the above bound as illustrated in Figure 3.

Our second contribution is a sharp upper bound of $s$ under which a simple Wald-type testing method based on $\bar{f}$ is minimax optimal in the sense of [6]. It is not surprising that our testing method is consistent no matter $s$ is fixed or diverges at any rate. Rather, this sharp bound is entirely determined by analyzing its (non-asymptotic) power. Specifically, we find that our testing method is minimax optimal if and only if $s$ does not grow faster than $N^{(4m-1)/(4m+1)}$. Again, we observe a subtle interplay between $s$ and $\lambda$ as depicted in Figure 3.

One theoretical insight obtained in our setup is that a more smooth regression function can be optimally estimated or tested at a shorter time. In addition, the above Figure 3 implies that $s$ and $\lambda$ play an interchangeable role in obtaining statistical optimality. Therefore, we argue that it might be attempting to view sample splitting as an alternative form of regularization, complementing the use of penalization in smoothing spline. In practice, we propose to select $\lambda$ via a distributed version of generalized cross validation (GCV); see [14].

In the end, we want to mention that our theoretical results are developed in one-dimensional models under fixed design. This setting allows us to develop proofs based on exact analysis of various Fourier series, coupled with properties of circulant Bernoulli polynomial kernel matrix. The major goal of this work is to provide some theoretical insights in a relatively simple setup, which are useful in extending our results to more general setup such as random or multi-dimensional design. Efforts toward this direction have been made by [8] who derived upper bounds of $s$ for optimal estimation or testing in various nonparametric models when design is random and multi-dimensional.

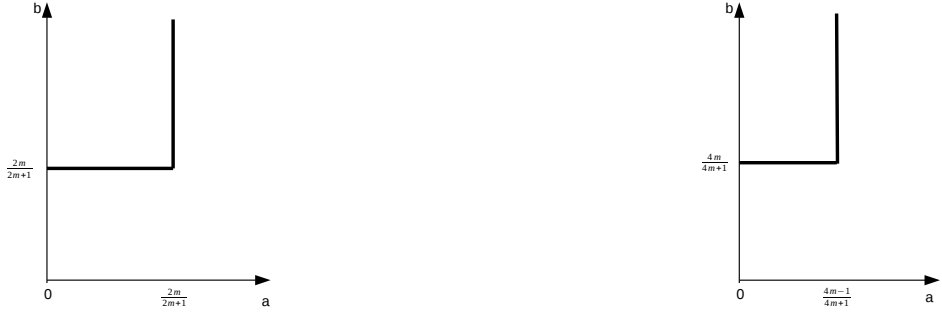

FIG 3. *Two lines indicate the choices of $s \asymp N^a$ and $\lambda \asymp N^{-b}$, leading to minimax optimal estimation rate (left) and minimax optimal testing rate (right). Whereas $(a, b)$'s outside these two lines lead to suboptimal rates. Results are based on smoothing spline regression with regularity $m \geq 1$.*

## 2. Smoothing Spline Model.

Suppose that we observe samples from model (1.1). The regression function $f$ is smooth in the sense that it belongs to an $m$-order ($m \geq 1$) periodic Sobolev space:

$$S^m(\mathbb{I}) = \left\{ \sum_{\nu=1}^{\infty} f_\nu \varphi_\nu(\cdot) : \sum_{\nu=1}^{\infty} f_\nu^2 \gamma_\nu < \infty \right\},$$

where $\mathbb{I} := [0, 1]$ and for $k = 1, 2, \ldots,$

$$\varphi_{2k-1}(t) = \sqrt{2} \cos(2\pi k t), \quad \varphi_{2k}(t) = \sqrt{2} \sin(2\pi k t),$$

$$\gamma_{2k-1} = \gamma_{2k} = (2\pi k)^{2m}.$$

The entire dataset is distributed to each machine in a uniform manner as follows. For $j = 1, \ldots, s$, the $j$th machine is assigned with samples $(Y_{i,j}, t_{i,j})$, where

$$Y_{i,j} = y_{is-s+j-1} \text{ and } t_{i,j} = \frac{is - s + j - 1}{N}$$

for $i = 1, \ldots, n$. Obviously, $t_{1,j}, \ldots, t_{n,j}$ are evenly spaced points (with a gap $1/n$) across $\mathbb{I}$. At the $j$th machine, we have the following sub-model:

$$(2.1) \qquad Y_{i,j} = f(t_{i,j}) + \epsilon_{i,j}, \ i = 1, \ldots, n,$$

where $\epsilon_{i,j} = \epsilon_{is-s+j-1}$, and obtain the $j$th sub-estimate as

$$\widehat{f}_j = \arg \min_{f \in S^m(\mathbb{I})} \ell_{j,n,\lambda}(f).$$

Here, $\ell_{j,n,\lambda}$ represents a penalized square criterion function based on the jth subsample:

$$(2.2) \qquad \ell_{j,n,\lambda}(f) = \frac{1}{2n} \sum_{i=1}^{n} (Y_{i,j} - f(t_{i,j}))^2 + \frac{\lambda}{2} J(f,f),$$

with $\lambda > 0$ being a smoothing parameter and $J(f,g) = \int_{\mathbb{I}} f^{(m)}(t) g^{(m)}(t) dt$[1]

## 3. Minimax Optimal Estimation.

In this section, we investigate the impact of the number of machines on the mean squared error of $\bar{f}$. Specifically, Theorem 3.1 provides an (non-asymptotic) upper bound for this mean squared error, while Theorem 3.2 provides a (non-asymptotic) lower bound. Notably, both bounds hold uniformly as $s$ diverges. From these bounds, we observe an interesting phase transition phenomenon that $\bar{f}$ is minimax optimal if $s$ does not grow faster than $N^{2m/(2m+1)}$ and an optimal $\lambda \asymp N^{-2m/(2m+1)}$ is chosen, but the minimax optimality breaks down if $s$ grows even slightly faster (no matter how $\lambda$ is chosen). Hence, the upper bound of $s$ is sharp. Moreover, $\lambda$ does not need to be optimal when this bound is attained. In some sense, a proper sample splitting can compensate a sub-optimal choice of $\lambda$.

In this section, we assume that $\epsilon_l$'s are *iid* zero-mean random variables with unit variance. Denote mean squared error as

$$\text{MSE}_{f_0}(f) := E_{f_0}\{\|f - f_0\|_2^2\},$$

where $\|f\|_2 = \sqrt{\int_{\mathbb{I}} f(t)^2 dt}$. For simplicity, we write $E_{f_0}$ as $E$ later. Define $h = \lambda^{1/(2m)}$.

THEOREM 3.1. *(Upper Bounds of Variance and Squared Bias) Suppose $h > 0$, and $N$ is divisible by $n$. Then there exist absolute positive constants $b_m, c_m \geq 1$ (depending on $m$ only) such that*

$$(3.1) \qquad E\{\|\bar{f} - E\{\bar{f}\}\|_2^2\} \leq b_m \left( N^{-1} + (Nh)^{-1} \int_0^{\pi nh} \frac{1}{(1+x^{2m})^2} dx \right),$$

$$(3.2) \qquad \|E\{\bar{f}\} - f_0\|_2 \leq c_m \sqrt{J(f_0)(\lambda + n^{-2m} + N^{-1})}$$

*for any fixed $1 \leq s \leq N$.*

From (6.2) and (6.3) in Appendix, we can tell that $\bar{f} - E\{\bar{f}\}$ is irrelevant to $f_0$. So is the upper bound for the (integrated) variance in (3.1). However, this is not the case for the (integrated) bias $\|E\{\bar{f}\} - f_0\|_2$, whose upper bound depends on $f_0$ through its norm $J(f_0)$. In particular, the (integrated) bias becomes zero if $f_0$ is in the null space, i.e., $J(f_0) = 0$, according to (3.2).

Since

$$(3.3) \qquad \text{MSE}_{f_0}(\bar{f}) = E\{\|\bar{f} - E\{\bar{f}\}\|_2^2\} + \|E\{\bar{f}\} - f_0\|_2^2,$$

Theorem 3.1 says that

$$(3.4) \qquad \text{MSE}_{f_0}(\bar{f}) \leq b_m \left( N^{-1} + (Nh)^{-1} \int_0^{\pi nh} \frac{1}{(1+x^{2m})^2} dx \right) + c_m^2 J(f_0)(\lambda + n^{-2m} + N^{-1}).$$

---

[1] For simplicity, we denote $J(f,f) = J(f)$ later.

When we choose $h \asymp N^{-1/(2m+1)}$ and $n^{-2m} = O(\lambda)$, it can be seen from (3.4) that $\bar{f}$ is minimax optimal, i.e., $\|\bar{f} - f_0\|_2 = O_P(N^{-m/(2m+1)})$. Obviously, the above two conditions hold if

$$(3.5) \qquad \lambda \asymp N^{-2m/(2m+1)} \ \text{ and } \ s = O(N^{2m/(2m+1)}).$$

From now on, we define the optimal choice of $\lambda$ as $N^{-2m/(2m+1)}$, denoted as $\lambda_*$; according to [16]. Alternatively, the minimax optimality can be achieved if $s \asymp N^{2m/(2m+1)}$ and $nh = o(1)$, i.e., $\lambda = o(\lambda_*)$. In other words, a sub-optimal choice of $\lambda$ can be compensated by a proper sampling splitting strategy. See Figure 3 for the subtle relation between $s$ and $\lambda$. It should be mentioned that $\lambda_*$ depends on $N$ (rather than $n$) for achieving optimal estimation rate. In practice, we propose to select $\lambda$ via a distributed version of GCV; see [14].

REMARK 3.1. *Under random design and uniformly bounded eigenfunctions, Corollary 4 in [16] showed that the above rate optimality is achieved under the following upper bound on $s$ (and $\lambda = \lambda_*$)*

$$s = O(N^{(2m-1)/(2m+1)}/\log N).$$

*For example, when $m = 2$, their upper bound is $N^{0.6}/\log N$ (versus $N^{0.8}$ in our case). We improve their upper bound by applying a more direct proof strategy.*

To understand whether our upper bound can be further improved, we prove a lower bound result in a "worst case" scenario. Specifically, Theorem 3.2 implies that once $s$ is beyond the above upper bound, the rate optimality will break down for at least one true $f_0$.

THEOREM 3.2. *(Lower Bound of Squared Bias) Suppose $h > 0$, and $N$ is divisible by $n$. Then for any constant $C > 0$, it holds that*

$$\sup_{\substack{f_0 \in S^m(\mathbb{I}) \\ J(f_0) \le C}} \|E\{\bar{f}\} - f_0\|_2^2 \ge C(a_m n^{-2m} - 8N^{-1}),$$

*where $a_m \in (0, 1)$ is an absolute constant depending on $m$ only, for any fixed $1 < s < N$.*

It follows by (3.3) that

$$(3.6) \qquad \sup_{\substack{f_0 \in S^m(\mathbb{I}) \\ J(f_0) \le C}} \mathrm{MSE}_{f_0}(\bar{f}) \ge \sup_{\substack{f_0 \in S^m(\mathbb{I}) \\ J(f_0) \le C}} \|E\{\bar{f}\} - f_0\|_2^2 \ge C(a_m n^{-2m} - 8N^{-1}).$$

It is easy to check that the above lower bound is strictly slower than the optimal rate $N^{-2m/(2m+1)}$ if $s$ grows faster than $N^{2m/(2m+1)}$ no matter how $\lambda$ is chosen. Therefore, we claim that $N^{2m/(2m+1)}$ is a sharp upper bound of $s$ for obtaining an averaged smoothing spline estimate.

In the end, we provide a graphical interpretation for our sharp bound result. Let $s = N^a$ for $0 \le a \le 1$ and $\lambda = N^{-b}$ for $0 < b < 2m$. Define $\rho_1(a), \rho_2(a)$ and $\rho_3(a)$ as

$$\text{Upper bound of squared bias: } N^{-\rho_1(a)} \;\asymp\; \lambda + n^{-2m} + N^{-1},$$

$$\text{Lower bound of squared bias: } N^{-\rho_2(a)} \;\asymp\; \max\{n^{-2m} - N^{-1}, 0\},$$

$$\text{Upper bound of variance: } N^{-\rho_3(a)} \;\asymp\; N^{-1} + (Nh)^{-1} \int_0^{\pi n h} \frac{1}{(1 + x^{2m})^2} \, dx,$$

based on Theorems 3.1 and 3.2. A direct examination reveals that

$$\rho_1(a) = \min\{2m(1-a), 1, b\}$$

$$\rho_2(a) = \begin{cases} 2m(1-a), & a > (2m-1)/(2m) \\ \infty, & a \le (2m-1)/(2m) \end{cases}$$

$$\rho_3(a) = \max\{a, (2m-b)/(2m)\}$$

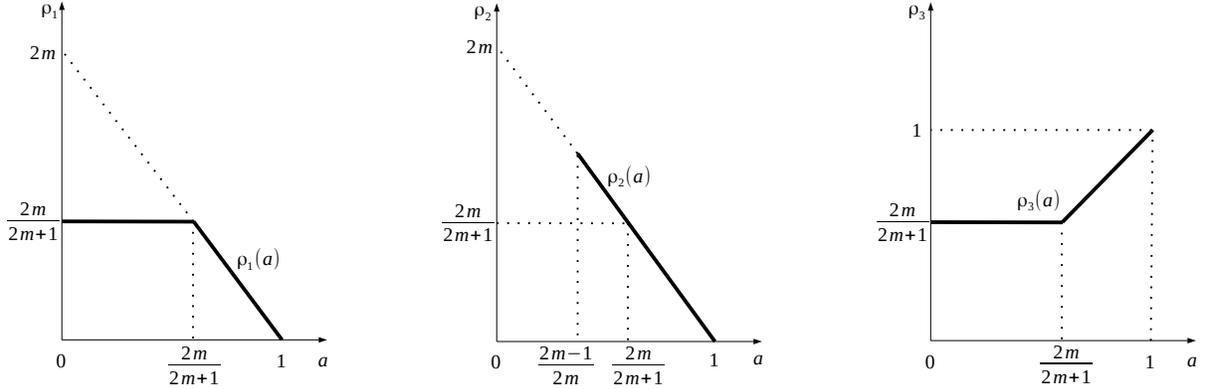

FIG 4. *Plots of $\rho_1(a), \rho_2(a), \rho_3(a)$ versus $a$, indicated by thick solid lines, under $\lambda = N^{-2m/(2m+1)}$. $\rho_1(a), \rho_2(a)$ and $\rho_3(a)$ indicate upper bound of squared bias, lower bound of squared bias and upper bound of variance, respectively. $\rho_2(a)$ is plotted only for $(2m-1)/(2m) < a \le 1$; when $0 \le a \le (2m-1)/(2m), \rho_2(a) = \infty$, which is omitted.*

Figure 4 displays $\rho_1, \rho_2, \rho_3$ for $\lambda = N^{-2m/(2m+1)}$. It can be seen that when $a \in [0, 2m/(2m+1)]$, upper bounds of squared bias and variance maintain at the same optimal rate $N^{-2m/(2m+1)}$, while the exact bound of squared bias increases above $N^{-2m/(2m+1)}$ when $a \in (2m/(2m+1), 1)$. This explains why transition occurs at the critical point $a = 2m/(2m+1)$ (even when the upper bound of variance decreases below $N^{-2m/(2m+1)}$ when $a \in (2m/(2m+1), 1)$).

It should be mentioned that when $\lambda \ne N^{-2m/(2m+1)}$, i.e., $b \ne 2m/(2m+1)$, suboptimal estimation almost always occurs. More explicitly, $b < 2m/(2m+1)$ yields $\rho_1(a) < 2m/(2m+1)$ for any $0 \le a \le 1$. While $b > 2m/(2m+1)$ yields $\rho_2(a) < 2m/(2m+1)$ for any $2m/(2m+1) < a \le 1$; yields $\rho_3(a) < 2m/(2m+1)$ for any $0 \le a < 2m/(2m+1)$. The only exception is $a = 2m/(2m+1)$ which yields $\rho_1 = \rho_2 = \rho_3 = 2m/(2m+1)$ for any $b > 2m/(2m+1)$.

REMARK 3.2. *As a side remark, we notice that each machine is assigned with $n \asymp N^{1/(2m+1)}$ samples when $s$ attains its upper bound in the estimation regime. This is very similar as the local* *polynomial estimation where approximately $N^{1/(2m+1)}$ local points are used for obtaining optimal estimation (although we realize that our data is distributed in a global manner).*

REMARK 3.3. *Under repeated curves with a common design, [2] observed a similar phase transition phenomenon for the minimax rate of a two-stage estimate, where the rate transits when the number of sample curves is nearly $N^{2m/(2m+1)}$. This coincides with our observation for $s$. However, the common design assumption, upon which their results crucially rely, clearly does not apply to our divide-and-conquer setup, and our proof techniques are significantly different. Rather, Theorems 3.1 and 3.2 imply that the results in [2] may still hold for a non-common design.*

## 4. Minimax Optimal Testing.
In this section, we consider nonparametric testing:

$$(4.1) \qquad H_0 : f = 0 \quad \text{v.s.} \quad H_1 : f \in S^m(\mathbb{I}) \backslash \{0\}.$$

In general, testing $f = f_0$ (for a known $f_0$) is equivalent to testing $f_* \equiv f - f_0 = 0$. So (4.1) has no loss of generality. Inspired by the classical Wald test ([11]), we propose a simple test statistic based on the $\bar{f}$ as

$$T_{N,\lambda} := \|\bar{f}\|_2^2.$$

We find that testing consistency essentially requires no condition on the number of machines no matter it is fixed or diverges at any rate. However, our power analysis, which is non-asymptotically valid, depends on the number of machines in a nontrivial way. Specifically, we discover that our test method is minimax optimal in the sense of Ingster ([6]) when $s$ does not grow faster than $N^{(4m-1)/(4m+1)}$ and $\lambda$ is chosen optimally (different from $\lambda_*$, though), but it is no longer optimal once $s$ is beyond the above threshold (no matter how $\lambda$ is chosen). This is a similar phase transition phenomenon as we observe in the estimation regime. Again, we notice an optimal choice of $\lambda$ may not be necessary if the above upper bound of $s$ is achieved.

In this section, we assume that the model errors $\epsilon_{i,j}$'s are *iid* standard normal for technical convenience. In fact, our results can be generalized to likelihood ratio test without assuming Gaussian errors. This extension is possible (technically tedious, though) since likelihood ratio statistic can be approximated by $T_{N,\lambda}$ through quadratic expansion; see [9].

Theorem 4.1 implies the consistency of our proposed test method with the following testing rule:

$$\phi_{N,\lambda} = I(|T_{N,\lambda} - \mu_{N,\lambda}| \geq z_{1-\alpha/2}\sigma_{N,\lambda}),$$

where $\mu_{N,\lambda} := \mathrm{E}_{H_0}\{T_{N,\lambda}\}$, $\sigma_{N,\lambda}^2 := \mathrm{Var}_{H_0}\{T_{N,\lambda}\}$ and $z_{1-\alpha/2}$ is the $(1 - \alpha/2) \times 100$ percentile of $N(0,1)$. The conditions required in Theorem 4.1 are so mild that our proposed testing is consistent no matter the number of machines is fixed or diverges at any rate.

THEOREM 4.1. *(Testing Consistency) Suppose that $h \to 0$, $n \to \infty$ when $N \to \infty$, and $\lim_{N \to \infty} nh$ exists (which could be infinity). Then, we have under $H_0$,*

$$\frac{T_{N,\lambda} - \mu_{N,\lambda}}{\sigma_{N,\lambda}} \xrightarrow{d} N(0,1), \quad as \ N \to \infty.$$

Our next theorem analyzes the non-asymptotic power of $T_{N,\lambda}$, in which we pay particular attention to the impact of $s$ on the separation rate of testing, defined as

$$d_{N,\lambda} = \sqrt{\lambda + n^{-2m} + \sigma_{N,\lambda}}.$$

Let $\mathcal{B} = \{f \in S^m(\mathbb{I}) : J(f) \le C\}$ for a positive constant $C$.

THEOREM 4.2. *(Upper Bound) Suppose that $h \to 0$, $n \to \infty$ when $N \to \infty$, and $\lim_{N \to \infty} nh$ exists (which could be infinity). Then for any $\varepsilon > 0$, there exist $C_\varepsilon, N_\varepsilon > 0$ s.t. for any $N \ge N_\varepsilon$,*

$$(4.2) \qquad \inf_{\substack{f \in \mathcal{B} \\ \|f\|_2 \ge C_\varepsilon d_{N,\lambda}}} P_f(\phi_{N,\lambda} = 1) \ge 1 - \varepsilon.$$

Under assumptions of Theorem 4.1, it can be shown that (see (6.40) in Appendix)

$$(4.3) \qquad \sigma^2_{N,\lambda} \asymp \begin{cases} \frac{n}{N^2}, & \text{if } \lim_{N \to 0} nh = 0, \\ \frac{1}{N^2 h}, & \text{if } \lim_{N \to \infty} nh > 0. \end{cases}$$

Given a range of $\lambda$ leading to $\lim_{N \to \infty} nh > 0$, we have by (4.3) that $d_{N,\lambda} = \sqrt{\lambda + (Nh^{1/2})^{-1}}$. An optimal choice of $\lambda$ (satisfying the above requirement) is $\lambda_{**} := N^{-4m/(4m+1)}$ since it leads to the optimal separating rate $d^*_{N,\lambda} := N^{-2m/(4m+1)}$; see [6]. Meanwhile, the constraint $\lim_{N \to \infty} nh > 0$ (together with the choice of $\lambda_{**}$) implies that

$$(4.4) \qquad s = O(N^{(4m-1)/(4m+1)}).$$

The above discussions illustrate that we can always choose $\lambda_{**}$ to obtain a minimax optimal testing (just as in the single dataset [9]) as long as $s$ does not grow faster than $N^{(4m-1)/(4m+1)}$. In the case that $\lim_{N \to \infty} nh = 0$, the minimax optimality can be maintained if $s \asymp N^{(4m-1)/(4m+1)}$, $h = o(1)$ and $nh = o(1)$. Such a selection of $s$ gives us a lot of freedom in choosing $\lambda$ that needs to satisfy $\lambda = o(\lambda_{**})$. A complete picture in depicting the relation between $s$ and $\lambda$ is given in Figure 3.

We further discover in Theorem 4.3 that the upper bound (4.4) turns out to be sharp.

THEOREM 4.3. *(Lower Bound) Suppose that $s \gg N^{(4m-1)/(4m+1)}$, $h \to 0$, $n \to \infty$ when $N \to \infty$, and $\lim_{N \to \infty} nh$ exists (which could be infinity). Then there exists a positive sequence $\beta_{N,\lambda}$ with $\lim_{N \to \infty} \beta_{N,\lambda} = \infty$ s.t.*

$$(4.5) \qquad \limsup_{N \to \infty} \inf_{\substack{f \in \mathcal{B} \\ \|f\|_2 \ge \beta_{N,\lambda} d^*_{N,\lambda}}} P_f(\phi_{N,\lambda} = 1) \le \alpha.$$

*Recall that $1 - \alpha$ is the pre-specified significance level.*

Theorem $4.3$ says that when $s \gg N^{(4m-1)/(4m+1)}$, the test $\phi_{N,\lambda}$ is no longer powerful even when $\|f\|_2 \gg d^*_{N,\lambda}$. In other words, our test method fails to be optimal. Therefore, we claim that $N^{(4m-1)/(4m+1)}$ is a sharp upper bound of $s$ to ensure our testing to be minimax optimal.

REMARK 4.1.   *As a side remark, the existence of* $\lim_{N\to\infty} nh$ *can be replaced by the following weaker condition under which the results in Theorems $4.1$, $4.2$ and $4.3$ still hold:*

$$Condition \ (\mathbf{R}): either \ \lim_{N\to\infty} \ nh = 0 \ or \ \inf_{N\geq 1} nh > 0.$$

*Condition ($\mathbf{R}$) aims to exclude irregularly behaved $s$ such as in the following case where $s$ vibrates too much along with $N$:*

$$(4.6) \qquad\qquad s = \begin{cases} N^{b_1}, & N \ is \ odd, \\ N^{b_2}, & N \ is \ even, \end{cases}$$

*where $h \asymp N^{-c}$ for some $c > 0$, $b_1, b_2 \in [0,1]$ satisfy $b_1 + c \geq 1$ and $b_2 + c < 1$. Clearly, Condition ($\mathbf{R}$) fails under ($4.6$).*

**5. Discussions.**  This paper offers "theoretical" suggestions on the allocation of data. In a relatively simple distributed algorithm, i.e., in $m$-order periodic splines with evenly spaced design, our recommendation proceeds as follows:

- Distribute to
$$s \asymp N^{2m/(2m+1)}$$
  machines for obtaining an optimal estimate;
- Distribute to
$$s \asymp N^{(4m-1)/(4m+1)}$$
  machines for performing an optimal test.

However, data-dependent formulae are still needed in picking a right number of machines in practice. This might be possible in light of Figure 3 indicating that sample splitting could be an alternative form of tuning. As for the choice of $\lambda$, we prove that it should be chosen in the order of $N$ even when each subsample has size $n$. Hence, a distributed version of the generalized cross validation method is applied to each sub-sample; see [14]. Another theoretically interesting direction is how much adaptive estimation (where $m$ is unknown) can affect the computational limits.

**Acknowledgments** We thank PhD student Meimei Liu at Purdue for the simulation study.

**6. Appendix.**  Proofs of our results are included in this section.

6.1. *Proofs in Section 3.*

PROOF OF THEOREM 3.1. We do a bit preliminary analysis before proving (3.1) and (3.2). It follows from [13] that $(S^m(\mathbb{I}), J)$ is a reproducing kernel Hilbert space with reproducing kernel function

$$K(x, y) = \sum_{\nu=1}^{\infty} \frac{\varphi_\nu(x)\varphi_\nu(y)}{\gamma_\nu} = 2\sum_{k=1}^{\infty} \frac{\cos(2\pi k(x - y))}{(2\pi k)^{2m}}, \ x, y \in \mathbb{I}.$$

For convenience, define $K_x(\cdot) = K(x, \cdot)$ for any $x \in \mathbb{I}$. It follows from the representer theorem ([13]) that the optimization to problem (2.2) has a solution

$$\widehat{f}_j = \sum_{i=1}^{n} \widehat{c}_{i,j} K_{t_{i,j}}, \ j = 1, 2, \ldots, s, \tag{6.1}$$

where $\widehat{c}_j = (\widehat{c}_{1,j}, \ldots, \widehat{c}_{n,j})^T = n^{-1}(\Sigma_j + \lambda I_n)^{-1}Y_j$, $Y_j = (Y_{1,j}, \ldots, Y_{n,j})^T$, $I_n$ is $n \times n$ identity matrix, and $\Sigma_j = [K(t_{i,j}, t_{i',j})/n]_{1 \le i, i' \le n}$. It is easy to see that $\Sigma_1 = \Sigma_2 = \cdots = \Sigma_s$. For convenience, denote $\Sigma = \Sigma_1$. Similarly, define

$$K'(x, y) = \sum_{\nu=1}^{\infty} \frac{\varphi_\nu(x)\varphi_\nu(y)}{\gamma_\nu^2} = 2\sum_{k=1}^{\infty} \frac{\cos(2\pi k(x - y))}{(2\pi k)^{4m}}, \ x, y \in \mathbb{I}.$$

For $1 \le j \le s$, let $\Omega_j = [K'(t_{i,j}, t_{i',j})/n]_{1 \le i, i' \le n}$. It is easy to see that $\Omega_1 = \Omega_2 = \cdots = \Omega_s$. For convenience, denote $\Omega = \Omega_1$, and let $\Phi_{\nu,j} = (\varphi_\nu(t_{1,j}), \ldots, \varphi_\nu(t_{n,j}))$.

It is easy to examine that

$$\begin{aligned}
\bar{f} &= \sum_{\nu=1}^{\infty} \frac{\sum_{j=1}^{s} \Phi_{\nu,j}(\Sigma + \lambda I_n)^{-1} Y_j}{N\gamma_\nu} \varphi_\nu \\
&= \sum_{\nu=1}^{\infty} \frac{\sum_{j=1}^{s} \Phi_{\nu,j}(\Sigma + \lambda I_n)^{-1}(\mathbf{f}_{0,j} + \boldsymbol{\epsilon}_j)}{N\gamma_\nu} \varphi_\nu,
\end{aligned} \tag{6.2}$$

and

$$E\{\bar{f}\} = \sum_{\nu=1}^{\infty} \frac{\sum_{j=1}^{s} \Phi_{\nu,j}(\Sigma + \lambda I_n)^{-1}\mathbf{f}_{0,j}}{N\gamma_\nu} \varphi_\nu, \tag{6.3}$$

where $\mathbf{f}_{0,j} = (f_0(t_{1,j}), \ldots, f_0(t_{n,j}))^T$ and $\boldsymbol{\epsilon}_j = (\epsilon_{1,j}, \ldots, \epsilon_{n,j})^T$.

We now look at $\Sigma$ and $\Omega$. For $0 \le l \le n - 1$, let

$$\begin{aligned}
c_l &= \frac{2}{n} \sum_{k=1}^{\infty} \frac{\cos(2\pi k l/n)}{(2\pi k)^{2m}}, \\
d_l &= \frac{2}{n} \sum_{k=1}^{\infty} \frac{\cos(2\pi k l/n)}{(2\pi k)^{4m}}.
\end{aligned}$$

Since $c_l = c_{n-l}$ and $d_l = d_{n-l}$ for $l = 1, 2, \ldots, n - 1$, $\Sigma$ and $\Omega$ are both symmetric circulant of order $n$. Let $\varepsilon = \exp(2\pi\sqrt{-1}/n)$. $\Omega$ and $\Sigma$ share the same normalized eigenvectors as

$$x_r = \frac{1}{\sqrt{n}}(1, \varepsilon^r, \varepsilon^{2r}, \ldots, \varepsilon^{(n-1)r})^T, \ r = 0, 1, \ldots, n - 1.$$

Let $M = (x_0, x_1, \ldots, x_{n-1})$. Denote $M^*$ as the conjugate transpose of $M$. Clearly, $MM^* = I_n$ and $\Sigma, \Omega$ admits the following decomposition

$$(6.4) \qquad \Sigma = M\Lambda_c M^*, \ \ \Omega = M\Lambda_d M^*,$$

where $\Lambda_c = \text{diag}(\lambda_{c,0}, \lambda_{c,1}, \ldots, \lambda_{c,n-1})$ and $\Lambda_d = \text{diag}(\lambda_{d,0}, \lambda_{d,1}, \ldots, \lambda_{d,n-1})$ with $\lambda_{c,l} = c_0 + c_1\varepsilon^l + \ldots + c_{n-1}\varepsilon^{(n-1)l}$ and $\lambda_{d,l} = d_0 + d_1\varepsilon^l + \ldots + d_{n-1}\varepsilon^{(n-1)l}$.

Direct calculations show that

$$(6.5) \qquad \lambda_{c,l} = \begin{cases} 2\sum_{k=1}^{\infty} \frac{1}{(2\pi kn)^{2m}}, & l = 0, \\ \sum_{k=1}^{\infty} \frac{1}{[2\pi(kn-l)]^{2m}} + \sum_{k=0}^{\infty} \frac{1}{[2\pi(kn+l)]^{2m}}, & 1 \le l \le n-1. \end{cases}$$

$$(6.6) \qquad \lambda_{d,l} = \begin{cases} 2\sum_{k=1}^{\infty} \frac{1}{(2\pi kn)^{4m}}, & l = 0, \\ \sum_{k=1}^{\infty} \frac{1}{[2\pi(kn-l)]^{4m}} + \sum_{k=0}^{\infty} \frac{1}{[2\pi(kn+l)]^{4m}}, & 1 \le l \le n-1. \end{cases}$$

It is easy to examine that

$$(6.7) \qquad \lambda_{c,0} = 2\bar{c}_m(2\pi n)^{-2m}, \ \ \lambda_{d,0} = 2\bar{d}_m(2\pi n)^{-4m},$$

and for $1 \le l \le n-1$,

$$\begin{aligned} \lambda_{c,l} &= \frac{1}{[2\pi(n-l)]^{2m}} + \frac{1}{(2\pi l)^{2m}} \\ &\quad + \sum_{k=2}^{\infty} \frac{1}{[2\pi(kn-l)]^{2m}} + \sum_{k=1}^{\infty} \frac{1}{[2\pi(kn+l)]^{2m}}, \\ \lambda_{d,l} &= \frac{1}{[2\pi(n-l)]^{4m}} + \frac{1}{(2\pi l)^{4m}} \\ &\quad + \sum_{k=2}^{\infty} \frac{1}{[2\pi(kn-l)]^{4m}} + \sum_{k=1}^{\infty} \frac{1}{[2\pi(kn+l)]^{4m}}, \end{aligned} \tag{6.8}$$

and for $\bar{c}_m := \sum_{k=1}^{\infty} k^{-2m}$, $\underline{c}_m := \sum_{k=2}^{\infty} k^{-2m}$, $\bar{d}_m := \sum_{k=1}^{\infty} k^{-4m}$, $\underline{d}_m := \sum_{k=2}^{\infty} k^{-4m}$,

$$\underline{c}_m(2\pi n)^{-2m} \le \sum_{k=2}^{\infty} \frac{1}{[2\pi(kn-l)]^{2m}} \le \bar{c}_m(2\pi n)^{-2m},$$

$$\underline{c}_m(2\pi n)^{-2m} \le \sum_{k=1}^{\infty} \frac{1}{[2\pi(kn+l)]^{2m}} \le \bar{c}_m(2\pi n)^{-2m},$$

$$\underline{d}_m(2\pi n)^{-4m} \le \sum_{k=2}^{\infty} \frac{1}{[2\pi(kn-l)]^{4m}} \le \bar{d}_m(2\pi n)^{-4m},$$

$$\underline{d}_m(2\pi n)^{-4m} \le \sum_{k=1}^{\infty} \frac{1}{[2\pi(kn+l)]^{4m}} \le \bar{d}_m(2\pi n)^{-4m}.$$

For simplicity, we denote $I = E\{\|\bar{f} - E\{\bar{f}\}\|_2^2\}$ and $II = \|E\{\bar{f}\} - f_0\|_2^2$. Hence, $\text{MSE}_{f_0}(\bar{f}) = I + II$.

**Proof of (3.1)**

Using (6.4) – (6.8), we get that

$$
\begin{aligned}
I &= \sum_{\nu=1}^{\infty} \frac{\sum_{j=1}^{s} E\{|\Phi_{\nu,j}(\Sigma + \lambda I_n)^{-1} \epsilon_j|^2\}}{N^2 \gamma_\nu^2} \\
&= \sum_{\nu=1}^{\infty} \frac{\sum_{j=1}^{s} \operatorname{trace}((\Sigma + \lambda I_n)^{-1} \Phi_{\nu,j}^T \Phi_{\nu,j}(\Sigma + \lambda I_n)^{-1})}{N^2 \gamma_\nu^2} \\
&= \frac{n}{N^2} \sum_{j=1}^{s} \operatorname{trace}\left((\Sigma + \lambda I_n)^{-1} \sum_{\nu=1}^{\infty} \frac{\Phi_{\nu,j}^T \Phi_{\nu,j}/n}{\gamma_\nu^2} (\Sigma + \lambda I_n)^{-1}\right) \\
&= \frac{n}{N^2} \sum_{j=1}^{s} \operatorname{trace}\left((\Sigma + \lambda I_n)^{-1} \Omega (\Sigma + \lambda I_n)^{-1}\right) \\
&= \frac{1}{N} \operatorname{trace}\left(M(\Lambda_c + \lambda I_n)^{-1} \Lambda_d (\Lambda_c + \lambda I_n)^{-1} M^*\right) \\
&= \frac{1}{N} \sum_{l=0}^{n-1} \frac{\lambda_{d,l}}{(\lambda + \lambda_{c,l})^2} \\
&\leq \frac{2\bar{d}_m}{N(2\bar{c}_m + (2\pi n)^{2m} \lambda)^2} \\
&\quad + (1 + \bar{d}_m) N^{-1} \sum_{l=1}^{n-1} \frac{(2\pi(n-l))^{-4m} + (2\pi l)^{-4m}}{(\lambda + (2\pi(n-l))^{-2m} + (2\pi l)^{-2m})^2} \\
&\leq \frac{2\bar{d}_m}{N(2\bar{c}_m + (2\pi n)^{2m} \lambda)^2} \\
&\quad + 2(1 + \bar{d}_m) N^{-1} \sum_{1 \leq l \leq n/2} \frac{(2\pi l)^{-4m} + (2\pi(n-l))^{-4m}}{(\lambda + (2\pi l)^{-2m} + (2\pi(n-l))^{-2m})^2} \\
&\leq \frac{2\bar{d}_m}{N(2\bar{c}_m + (2\pi n)^{2m} \lambda)^2} + 4(1 + \bar{d}_m) N^{-1} \sum_{1 \leq l \leq n/2} \frac{(2\pi l)^{-4m}}{(\lambda + (2\pi l)^{-2m})^2} \\
&\leq \frac{2\bar{d}_m}{N(2\bar{c}_m + (2\pi n)^{2m} \lambda)^2} + \frac{2(1 + \bar{d}_m)}{\pi Nh} \int_0^{\pi nh} \frac{1}{(1 + x^{2m})^2} dx \\
&\leq b_m \left(\frac{1}{N} + \frac{1}{Nh} \int_0^{\pi nh} \frac{1}{(1 + x^{2m})^2} dx\right),
\end{aligned}
$$

where $b_m \geq 1$ is an absolute constant depending on $m$ only. This proves (3.1).

**Proof of (3.2)**

Throughout, let $\eta = \exp(2\pi\sqrt{-1}/N)$. For $1 \leq j, l \leq s$, define

$$
\begin{aligned}
\Sigma_{j,l} &= \frac{1}{n} \sum_{\nu=1}^{\infty} \frac{\Phi_{\nu,j}^T \Phi_{\nu,l}}{\gamma_\nu}, \\
\sigma_{j,l,r} &= \frac{2}{n} \sum_{k=1}^{\infty} \frac{\cos\left(2\pi k \left(\frac{r}{n} - \frac{j-l}{N}\right)\right)}{(2\pi k)^{2m}}, \ r = 0, 1, \ldots, n-1.
\end{aligned}
$$

It can be shown that $\Sigma_{j,l}$ is a circulant matrix with elements $\sigma_{j,l,0}, \sigma_{j,l,1}, \ldots, \sigma_{j,l,n-1}$, therefore, by [1] we get that

$$\Sigma_{j,l} = M\Lambda_{j,l}M^*, \tag{6.9}$$

where $M$ is the same as in (6.4), and $\Lambda_{j,l} = \mathrm{diag}(\lambda_{j,l,0}, \lambda_{j,l,1}, \ldots, \lambda_{j,l,n-1})$, with $\lambda_{j,l,r}$, for $r = 1, \ldots, n-1$, given by the following

$$
\begin{aligned}
\lambda_{j,l,r} &= \sum_{t=0}^{n-1} \sigma_{j,l,t}\varepsilon^{rt} \\
&= \frac{2}{n}\sum_{t=0}^{n-1}\sum_{k=1}^{\infty}\frac{\cos\left(2\pi k\left(\frac{t}{n} - \frac{j-l}{N}\right)\right)}{(2\pi k)^{2m}}\varepsilon^{rt} \\
&= \frac{1}{n}\sum_{k=1}^{\infty}\frac{\eta^{-k(j-l)}\sum_{t=0}^{n-1}\varepsilon^{(k+r)t} + \eta^{k(j-l)}\sum_{t=0}^{n-1}\varepsilon^{(r-k)t}}{(2\pi k)^{2m}} \\
&= \sum_{q=1}^{\infty}\frac{\eta^{-(qn-r)(j-l)}}{[2\pi(qn-r)]^{2m}} + \sum_{q=0}^{\infty}\frac{\eta^{(qn+r)(j-l)}}{[2\pi(qn+r)]^{2m}},
\end{aligned}
\tag{6.10}
$$

and for $r = 0$, given by

$$
\begin{aligned}
\lambda_{j,l,0} &= \sum_{t=0}^{n-1} \sigma_{j,l,t} \\
&= \frac{1}{n}\sum_{k=1}^{\infty}\frac{\sum_{t=0}^{n-1}\varepsilon^{kt}\eta^{k(j-l)} + \sum_{t=0}^{n-1}\varepsilon^{-kt}\eta^{-k(j-l)}}{(2\pi k)^{2m}} \\
&= \sum_{q=1}^{\infty}\frac{\eta^{qn(j-l)} + \eta^{-qn(j-l)}}{(2\pi qn)^{2m}}.
\end{aligned}
\tag{6.11}
$$

For $p \geq 0$, $1 \leq v \leq n$, $0 \leq r \leq n-1$ and $1 \leq j \leq s$, define

$$A_{p,v,r,j} = \frac{1}{s}\sum_{l=1}^{s}\lambda_{j,l,r}x_r^*\Phi_{2(pn+v)-1,l}^T, \quad B_{p,v,r,j} = \frac{1}{s}\sum_{l=1}^{s}\lambda_{j,l,r}x_r^*\Phi_{2(pn+v),l}^T.$$

By direct calculation, we have for $1 \leq v \leq n-1$,

$$
\begin{aligned}
\Phi_{2(pn+v)-1,l}x_r &= \sqrt{n/2}\left(\eta^{(pn+v)(l-1)}I(r+v=n) + \eta^{-(pn+v)(l-1)}I(v=r)\right), \\
\Phi_{2(pn+v),l}x_r &= \sqrt{-n/2}\left(\eta^{(pn+v)(l-1)}I(r+v=n) - \eta^{-(pn+v)(l-1)}I(v=r)\right),
\end{aligned}
\tag{6.12}
$$

and

$$
\begin{aligned}
\Phi_{2(pn+n)-1,l}x_r &= \sqrt{n/2}I(r=0)\left(\eta^{(p+1)n(l-1)} + \eta^{-(p+1)n(l-1)}\right), \\
\Phi_{2(pn+n),l}x_r &= \sqrt{-n/2}I(r=0)\left(\eta^{(p+1)n(l-1)} - \eta^{-(p+1)n(l-1)}\right).
\end{aligned}
\tag{6.13}
$$

Let $I(\cdot)$ be an indicator function. Then we have for $p \geq 0$, $1 \leq j \leq s$ and $1 \leq v, r \leq n-1$,

$$
\begin{aligned}
B_{p,v,r,j} &= \frac{1}{s} \sum_{l=1}^{s} \lambda_{j,l,r} x_r^* \Phi_{2(pn+v),l}^T \\
&= -\frac{1}{s} \sqrt{-n/2} \sum_{l=1}^{s} \left( \sum_{q=1}^{\infty} \frac{\eta^{-(qn-r)(j-l)}}{[2\pi(qn-r)]^{2m}} + \sum_{q=0}^{\infty} \frac{\eta^{(qn+r)(j-l)}}{[2\pi(qn+r)]^{2m}} \right) \\
&\quad \times \left( \eta^{-(pn+v)(l-1)} I(r+v=n) - \eta^{(pn+v)(l-1)} I(r=v) \right) \\
&= -\sqrt{-n/2} \left( \sum_{u \geq -p/s} \frac{\eta^{-(pn+v)(j-1)}}{[2\pi(uN+pn+v)]^{2m}} I(r+v=n) \right. \\
&\quad - \sum_{u \geq (p+1)/s} \frac{\eta^{(pn+v)(j-1)}}{[2\pi(uN-pn-v)]^{2m}} I(r=v) \\
&\quad + \sum_{u \geq (p+1)/s} \frac{\eta^{-(pn+v)(j-1)}}{[2\pi(uN-pn-v)]^{2m}} I(r+v=n) \\
&\quad \left. - \sum_{u \geq -p/s} \frac{\eta^{(pn+v)(j-1)}}{[2\pi(uN+pn+v)]^{2m}} I(r=v) \right) = a_{p,v} x_r^* \Phi_{2(pn+v),j}^T,
\end{aligned}
$$

(6.14)

where $a_{p,v} = \sum_{u \geq -p/s} \frac{1}{[2\pi(uN+pn+v)]^{2m}} + \sum_{u \geq (p+1)/s} \frac{1}{[2\pi(uN-pn-v)]^{2m}}$, for $p \geq 0$, $1 \leq v \leq n-1$.

For $v = n$, similar calculations give that

$$
\begin{aligned}
B_{p,n,r,j} &= -\sqrt{-n/2} I(r=0) \left( \sum_{u \geq -p/s} \frac{\eta^{-(pn+n)(j-1)}}{[2\pi(uN+pn+n)]^{2m}} \right. \\
&\quad - \sum_{u \geq (p+2)/s} \frac{\eta^{(pn+n)(j-1)}}{[2\pi(uN-pn-n)]^{2m}} \\
&\quad \left. + \sum_{u \geq (p+2)/s} \frac{\eta^{-(pn+n)(j-1)}}{[2\pi(uN-pn-n)]^{2m}} - \sum_{u \geq -p/s} \frac{\eta^{(pn+n)(j-1)}}{[2\pi(uN+pn+n)]^{2m}} \right) \\
&= a_{p,n} x_r^* \Phi_{2(pn+n),j}^T,
\end{aligned}
$$

(6.15)

where $a_{p,n} = \sum_{u \geq -p/s} \frac{1}{[2\pi(uN+pn+n)]^{2m}} + \sum_{u \geq (p+2)/s} \frac{1}{[2\pi(uN-pn-n)]^{2m}}$, for $p \geq 0$.

Similarly, we have $p \geq 0$, $1 \leq j \leq s$ and $1 \leq v, r \leq n-1$,

$$
\begin{aligned}
A_{p,v,r,j} &= \sqrt{n/2} \Bigg( \sum_{u \geq -p/s} \frac{\eta^{-(p n+v)(j-1)}}{[2\pi(uN+pn+v)]^{2m}} I(r+v=n) \\
&\quad + \sum_{u \geq (p+1)/s} \frac{\eta^{(pn+v)(j-1)}}{[2\pi(uN-pn-v)]^{2m}} I(r=v) \\
&\quad + \sum_{u \geq (p+1)/s} \frac{\eta^{-(pn+v)(j-1)}}{[2\pi(uN-pn-v)]^{2m}} I(r+v=n) \\
&\quad + \sum_{u \geq -p/s} \frac{\eta^{(pn+v)(j-1)}}{[2\pi(uN+pn+v)]^{2m}} I(r=v) \Bigg) \\
&= a_{p,v} x_r^* \Phi_{2(pn+v)-1,j}^T,
\end{aligned}
\tag{6.16}
$$

and for $v = n$,

$$
\begin{aligned}
A_{p,n,r,j} &= \sqrt{n/2}\, I(r=0) \Bigg( \sum_{u \geq -p/s} \frac{\eta^{-(pn+n)(j-1)}}{[2\pi(uN+pn+n)]^{2m}} \\
&\quad + \sum_{u \geq (p+2)/s} \frac{\eta^{(pn+n)(j-1)}}{[2\pi(uN-pn-n)]^{2m}} \\
&\quad + \sum_{u \geq (p+2)/s} \frac{\eta^{-(pn+n)(j-1)}}{[2\pi(uN-pn-n)]^{2m}} \\
&\quad + \sum_{u \geq -p/s} \frac{\eta^{(pn+n)(j-1)}}{[2\pi(uN+pn+n)]^{2m}} \Bigg) = a_{p,n} x_r^* \Phi_{2(pn+n)-1,j}^T.
\end{aligned}
\tag{6.17}
$$

It is easy to check that both (6.14) and (6.16) hold for $r = 0$. Summarizing (6.14)–(6.17), we have that for $p \geq 0$, $1 \leq j \leq s$, $1 \leq v \leq n$ and $0 \leq r \leq n-1$,

$$
\begin{aligned}
A_{p,v,r,j} &= a_{p,v} x_r^* \Phi_{2(pn+v)-1,j}^T, \\
B_{p,v,r,j} &= a_{p,v} x_r^* \Phi_{2(pn+v),j}^T.
\end{aligned}
\tag{6.18}
$$

To show (3.2), let $\bar{\boldsymbol{f}}_j = (E\{\bar{f}(t_{1,j})\}, \ldots, E\{\bar{f}(t_{n,j})\})^T$, for $1 \leq j \leq s$. It follows by (6.3) that

$$
\begin{aligned}
\bar{\boldsymbol{f}}_j &= \sum_{\nu=1}^{\infty} \frac{\sum_{l=1}^{s} \Phi_{\nu,l}(\Sigma + \lambda I_n)^{-1} \mathbf{f}_{0,l}}{N\gamma_\nu} \Phi_{\nu,j}^T \\
&= \frac{1}{s} \sum_{l=1}^{s} \left( \frac{1}{n} \sum_{\nu=1}^{\infty} \frac{\Phi_{\nu,j}^T \Phi_{\nu,l}}{\gamma_\nu} \right) (\Sigma + \lambda I_n)^{-1} \mathbf{f}_{0,l} \\
&= \frac{1}{s} \sum_{l=1}^{s} \Sigma_{j,l} (\Sigma + \lambda I_n)^{-1} \mathbf{f}_{0,l} \\
&= \frac{1}{s} \sum_{l=1}^{s} M \Lambda_{j,l} (\Lambda_c + \lambda I_n)^{-1} M^* \mathbf{f}_{0,l},
\end{aligned}
$$

together with (6.18), leading to that

$$
\begin{aligned}
M^* \bar{\boldsymbol{f}}_j &= \frac{1}{s} \sum_{l=1}^{s} \Lambda_{j,l} (\Lambda_c + \lambda I_n)^{-1} M^* \mathbf{f}_{0,l} \\
&= \sum_{\mu=1}^{\infty} f_\mu^0 \begin{pmatrix} \frac{\frac{1}{s} \sum_{l=1}^{s} \lambda_{j,l,0} x_0^* \Phi_{\mu,l}^T}{\lambda + \lambda_{c,0}} \\ \vdots \\ \frac{\frac{1}{s} \sum_{l=1}^{s} \lambda_{j,l,n-1} x_{n-1}^* \Phi_{\mu,l}^T}{\lambda + \lambda_{c,n-1}} \end{pmatrix} \\
&= \sum_{p=0}^{\infty} \sum_{v=1}^{n} f_{2(pn+v)-1}^0 \begin{pmatrix} \frac{\frac{1}{s} \sum_{l=1}^{s} \lambda_{j,l,0} x_0^* \Phi_{2(pn+v)-1,l}^T}{\lambda + \lambda_{c,0}} \\ \vdots \\ \frac{\frac{1}{s} \sum_{l=1}^{s} \lambda_{j,l,n-1} x_{n-1}^* \Phi_{2(pn+v)-1,l}^T}{\lambda + \lambda_{c,n-1}} \end{pmatrix} \\
&\quad + \sum_{p=0}^{\infty} \sum_{v=1}^{n} f_{2(pn+v)}^0 \begin{pmatrix} \frac{\frac{1}{s} \sum_{l=1}^{s} \lambda_{j,l,0} x_0^* \Phi_{2(pn+v),l}^T}{\lambda + \lambda_{c,0}} \\ \vdots \\ \frac{\frac{1}{s} \sum_{l=1}^{s} \lambda_{j,l,n-1} x_{n-1}^* \Phi_{2(pn+v),l}^T}{\lambda + \lambda_{c,n-1}} \end{pmatrix} \\
&= \sum_{p=0}^{\infty} \sum_{v=1}^{n} f_{2(pn+v)-1}^0 \begin{pmatrix} \frac{A_{p,v,0,j}}{\lambda + \lambda_{c,0}} \\ \vdots \\ \frac{A_{p,v,n-1,j}}{\lambda + \lambda_{c,n-1}} \end{pmatrix} \\
&\quad + \sum_{p=0}^{\infty} \sum_{v=1}^{n} f_{2(pn+v)}^0 \begin{pmatrix} \frac{B_{p,v,0,j}}{\lambda + \lambda_{c,0}} \\ \vdots \\ \frac{B_{p,v,n-1,j}}{\lambda + \lambda_{c,n-1}} \end{pmatrix} \\
&= \sum_{p=0}^{\infty} \sum_{v=1}^{n} f_{2(pn+v)-1}^0 \begin{pmatrix} \frac{a_{p,v}}{\lambda + \lambda_{c,0}} x_0^* \Phi_{2(pn+v)-1,j}^T \\ \vdots \\ \frac{a_{p,v}}{\lambda + \lambda_{c,n-1}} x_{n-1}^* \Phi_{2(pn+v)-1,j}^T \end{pmatrix} \\
&\quad + \sum_{p=0}^{\infty} \sum_{v=1}^{n} f_{2(pn+v)}^0 \begin{pmatrix} \frac{a_{p,v}}{\lambda + \lambda_{c,0}} x_0^* \Phi_{2(pn+v),j}^T \\ \vdots \\ \frac{a_{p,v}}{\lambda + \lambda_{c,n-1}} x_{n-1}^* \Phi_{2(pn+v),j}^T \end{pmatrix} .
\end{aligned}
$$

On the other hand,

$$
\begin{aligned}
M^* \mathbf{f}_{0,j} &= \sum_{\mu=1}^{\infty} f_\mu^0 M^* \Phi_{\mu,j}^T \\
&= \sum_{p=0}^{\infty} \sum_{v=1}^{n} f_{2(pn+v)-1}^0 M^* \Phi_{2(pn+v)-1,j}^T + \sum_{p=0}^{\infty} \sum_{v=1}^{n} f_{2(pn+v)}^0 M^* \Phi_{2(pn+v),j}^T \\
&= \sum_{p=0}^{\infty} \sum_{v=1}^{n} f_{2(pn+v)-1}^0 \left( \begin{array}{c} x_0^* \Phi_{2(pn+v)-1,j}^T \\ \vdots \\ x_{n-1}^* \Phi_{2(pn+v)-1,j}^T \end{array} \right) \\
&\quad + \sum_{p=0}^{\infty} \sum_{v=1}^{n} f_{2(pn+v)}^0 \left( \begin{array}{c} x_0^* \Phi_{2(pn+v),j}^T \\ \vdots \\ x_{n-1}^* \Phi_{2(pn+v),j}^T \end{array} \right).
\end{aligned}
$$

Therefore,

$$
\begin{aligned}
M^*(\bar{\boldsymbol{f}}_j - \mathbf{f}_{0,j}) &= \sum_{p=0}^{\infty} \sum_{v=1}^{n} f_{2(pn+v)-1}^0 \left( \begin{array}{c} b_{p,v,0} x_0^* \Phi_{2(pn+v)-1,j}^T \\ \vdots \\ b_{p,v,n-1} x_{n-1}^* \Phi_{2(pn+v)-1,j}^T \end{array} \right) \\
&\quad + \sum_{p=0}^{\infty} \sum_{v=1}^{n} f_{2(pn+v)}^0 \left( \begin{array}{c} b_{p,v,0} x_0^* \Phi_{2(pn+v),j}^T \\ \vdots \\ b_{p,v,n-1} x_{n-1}^* \Phi_{2(pn+v),j}^T \end{array} \right),
\end{aligned}
\tag{6.19}
$$

where $b_{p,v,r} = \frac{a_{p,v}}{\lambda + \lambda_{c,r}} - 1$, for $p \geq 0$, $1 \leq v \leq n$ and $0 \leq r \leq n-1$.

It holds the trivial observation $b_{ks+g,v,r} = b_{g,v,r}$ for $k \geq 0$, $0 \leq g \leq s-1$, $1 \leq v \leq n$ and $0 \leq r \leq n-1$. Define $C_{g,r} = \sum_{k=0}^{\infty} (f_{2(kN+gn+n-r)-1}^0 - \sqrt{-1} f_{2(kN+gn+n-r)}^0)$ and $D_{g,r} = \sum_{k=0}^{\infty} (f_{2(kN+gn+r)-1}^0 + \sqrt{-1} f_{2(kN+gn+r)}^0)$, for $0 \leq g \leq s-1$ and $0 \leq r \leq n-1$. Also denote $\overline{C_{g,r}}$ and $\overline{D_{g,r}}$ as their conjugate. By (6.12) and (6.13), and direct calculations we get that, for $1 \leq j \leq s$ and $1 \leq r \leq n-1$,

$$
\begin{aligned}
\delta_{j,r} &\equiv \sum_{p=0}^{\infty} \left( \sum_{v=1}^{n} f_{2(pn+v)-1}^0 b_{p,v,r} x_r^* \Phi_{2(pn+v)-1,j}^T + \sum_{v=1}^{n} f_{2(pn+v)}^0 b_{p,v,r} x_r^* \Phi_{2(pn+v),j}^T \right) \\
&= \sqrt{\frac{n}{2}} \sum_{p=0}^{\infty} \Big[ (f_{2(pn+n-r)-1}^0 - \sqrt{-1} f_{2(pn+n-r)}^0) b_{p,n-r,r} \eta^{-(pn+n-r)(j-1)} \\
&\quad + (f_{2(pn+r)-1}^0 + \sqrt{-1} f_{2(pn+r)}^0) b_{p,r,r} \eta^{(pn+r)(j-1)} \Big],
\end{aligned}
\tag{6.20}
$$

leading to that

$$
\begin{aligned}
\sum_{j=1}^{s} |\delta_{j,r}|^2 &= \frac{n}{2} \sum_{j=1}^{s} \left| \sum_{p=0}^{\infty} \Big[ (f^0_{2(pn+n-r)-1} - \sqrt{-1} f^0_{2(pn+n-r)}) b_{p,n-r,r} \eta^{-(pn+n-r)(j-1)} \right. \\
&\qquad \left. + (f^0_{2(pn+r)-1} + \sqrt{-1} f^0_{2(pn+r)}) b_{p,r,r} \eta^{(pn+r)(j-1)} \Big] \right|^2 \\
&= \frac{n}{2} \sum_{j=1}^{s} \left| \sum_{g=0}^{s-1} \Big( C_{g,r} b_{g,n-r,r} \eta^{-(gn+n-r)(j-1)} + D_{g,r} b_{g,r,r} \eta^{(gn+r)(j-1)} \Big) \right|^2 \\
&= \frac{n}{2} \sum_{g,g'=0}^{s-1} \sum_{j=1}^{s} (C_{g,r} b_{g,n-r,r} \eta^{-(gn+n-r)(j-1)} + D_{g,r} b_{g,r,r} \eta^{(gn+r)(j-1)}) \\
&\qquad \times (\overline{C_{g',r}} b_{g',n-r,r} \eta^{(g'n+n-r)(j-1)} + \overline{D_{g',r}} b_{g',r,r} \eta^{-(g'n+r)(j-1)}) \\
&= \frac{N}{2} \sum_{g=0}^{s-1} (|C_{g,r}|^2 b^2_{g,n-r,r} + C_{g,r} \overline{D_{s-1-g,r}} b_{g,n-r,r} b_{s-1-g,r,r} \\
&\qquad + D_{g,r} \overline{C_{s-1-g,r}} b_{g,r,r} b_{s-1-g,n-r,r} + |D_{g,r}|^2 b^2_{g,r,r}) \\
&= \frac{N}{2} \sum_{g=0}^{s-1} |C_{g,r} b_{g,n-r,r} + D_{s-1-g,r} b_{s-1-g,r,r}|^2 \\
&\le N \sum_{g=0}^{s-1} (|C_{g,r}|^2 b^2_{g,n-r,r} + |D_{s-1-g,r}|^2 b^2_{s-1-g,r,r}) \\
&= N \sum_{g=0}^{s-1} (|C_{g,r}|^2 b^2_{g,n-r,r} + |D_{g,r}|^2 b^2_{g,r,r}).
\end{aligned}
$$

(6.21)

It is easy to see that for $0 \le g \le s-1$ and $1 \le r \le n-1$,

$$
\begin{aligned}
|C_{g,r}|^2 &= (\sum_{k=0}^{\infty} f^0_{2(kN+gn+n-r)-1})^2 + (\sum_{k=0}^{\infty} f^0_{2(kN+gn+n-r)})^2 \\
&\le \sum_{k=0}^{\infty} (|f^0_{2(kN+gn+n-r)-1}|^2 + |f^0_{2(kN+gn+n-r)}|^2)(kN+gn+n-r)^{2m} \\
&\qquad \times \sum_{k=0}^{\infty} (kN+gn+n-r)^{-2m} \\
&\le \sum_{k=0}^{\infty} (|f^0_{2(kN+gn+n-r)-1}|^2 + |f^0_{2(kN+gn+n-r)}|^2)(kN+gn+n-r)^{2m} \\
&\qquad \times \frac{2m}{2m-1} (gn+n-r)^{-2m},
\end{aligned}
$$

(6.22)

and

$$
\begin{aligned}
|D_{g,r}|^2 &\le \sum_{k=0}^{\infty} (|f^0_{2(kN+gn+r)}|^2 + |f^0_{2(kN+gn+r)-1}|^2)(kN+gn+r)^{2m} \\
&\qquad \times \frac{2m}{2m-1} (gn+r)^{-2m}.
\end{aligned}
$$

(6.23)

For $1 \leq g \leq s-1$, we have $a_{g,n-r} \leq \lambda_{c,r}$, which further leads to $|b_{g,n-r,r}| \leq 2$. Meanwhile, by (6.5), we have

$$0 \leq \lambda_{c,r} - a_{0,r} \leq (2\pi(n-r))^{-2m} + 2\bar{c}_m(2\pi n)^{-2m} \leq (1 + 2\bar{c}_m)(2\pi(n-r))^{-2m}.$$

Then we have

$$
\begin{aligned}
|b_{0,r,r}| &= \frac{\lambda + \lambda_{c,r} - a_{0,r}}{\lambda + \lambda_{c,r}} \\
&\leq \frac{\lambda + (1 + 2\bar{c}_m)(2\pi(n-r))^{-2m}}{\lambda + (2\pi r)^{-2m} + (2\pi(n-r))^{-2m}} \\
&\leq (1 + 2\bar{c}_m)\frac{\lambda + (2\pi(n-r))^{-2m}}{\lambda + (2\pi r)^{-2m} + (2\pi(n-r))^{-2m}},
\end{aligned}
$$

leading to

$$
\begin{aligned}
r^{-2m}b_{0,r,r}^2 &\leq r^{-2m}(1 + 2\bar{c}_m)^2 \left( \frac{\lambda + (2\pi(n-r))^{-2m}}{\lambda + (2\pi(n-r))^{-2m} + (2\pi r)^{-2m}} \right)^2 \\
&\leq r^{-2m}(1 + 2\bar{c}_m)^2 \left( \frac{\lambda + (2\pi(n-r))^{-2m}}{\lambda + (2\pi(n-r))^{-2m} + (2\pi r)^{-2m}} \right) \\
&\leq (2\pi)^{2m}(1 + 2\bar{c}_m)^2(\lambda + (\pi n)^{-2m}).
\end{aligned}
\tag{6.24}
$$

The last inequality can be proved in two different cases: $2r \leq n$ and $2r > n$. Similarly, it can be shown that $(n-r)^{-2m}b_{0,n-r,r}^2 \leq (2\pi)^{2m}(1 + 2\bar{c}_m)^2(\lambda + (\pi n)^{-2m})$.

Then we have by (6.22)–(6.24) that

$$
\begin{aligned}
\sum_{r=1}^{n-1}\sum_{g=0}^{s-1}|D_{g,r}|^2 b_{g,r,r}^2 &\leq \sum_{r=1}^{n-1}\sum_{g=1}^{s-1}\sum_{k=0}^{\infty}(|f_{2(kN+gn+r)}^0|^2 + |f_{2(kN+gn+r)-1}^0|^2)(kN + gn + r)^{2m} \\
&\quad \times \frac{2m}{2m-1}(gn + r)^{-2m}2^2 + \sum_{r=1}^{n-1}\sum_{k=0}^{\infty}(|f_{2(kN+r)}^0|^2 + |f_{2(kN+r)-1}^0|^2)(kN + r)^{2m} \\
&\quad \times \frac{2m}{2m-1}r^{-2m}b_{0,r,r}^2 \\
&\leq c_m'(\lambda + n^{-2m})\sum_{r=1}^{n-1}\sum_{g=0}^{s-1}\sum_{k=0}^{\infty}(|f_{2(kN+gn+r)}^0|^2 + |f_{2(kN+gn+r)-1}^0|^2) \\
&\quad \times (2\pi(kN + gn + r))^{2m},
\end{aligned}
\tag{6.25}
$$

where $c_m' = \max\{(2\pi)^{-2m}\frac{8m}{2m-1}, (1 + 2\bar{c}_m)^2\frac{2m}{2m-1}\}$. Similarly, one can show that

$$
\begin{aligned}
\sum_{r=1}^{n-1}\sum_{g=0}^{s-1}|C_{g,r}|^2 b_{g,n-r,r}^2 &\leq c_m'(\lambda + n^{-2m})\sum_{r=1}^{n-1}\sum_{g=0}^{s-1}\sum_{k=0}^{\infty}(|f_{2(kN+gn+r)}^0|^2 + |f_{2(kN+gn+r)-1}^0|^2) \\
&\quad \times (2\pi(kN + gn + r))^{2m}.
\end{aligned}
\tag{6.26}
$$

Combining (6.25) and (6.26) we get that

$$
\sum_{r=1}^{n-1}\sum_{j=1}^{s}|\delta_{j,r}|^2 \;\leq\; 2c_m'(\lambda+n^{-2m})N\sum_{r=1}^{n-1}\sum_{g=0}^{s-1}\sum_{k=0}^{\infty}(|f^0_{2(kN+gn+r)}|^2+|f^0_{2(kN+gn+r)-1}|^2)
$$
$$
\times(2\pi(kN+gn+r))^{2m}. \tag{6.27}
$$

To the end of proof of (3.2), by (6.19) we have for $1\leq j\leq s$,

$$
\begin{aligned}
\delta_{j,0} &\equiv \sum_{p=0}^{\infty}\left(\sum_{v=1}^{n}f^0_{2(pn+v)-1}b_{p,v,0}x_0^*\Phi^T_{2(pn+v)-1,j}+\sum_{v=1}^{n}f^0_{2(pn+v)}b_{p,v,0}x_0^*\Phi^T_{2(pn+v),j}\right)\\
&= \sum_{p=0}^{\infty}\left(f^0_{2(pn+n)-1}b_{p,n,0}x_0^*\Phi^T_{2(pn+n)-1,j}+f^0_{2(pn+n)}b_{p,n,0}x_0^*\Phi^T_{2(pn+n),j}\right)\\
&= \sqrt{\frac{n}{2}}\sum_{p=0}^{\infty}\Big[(f^0_{2(pn+n)-1}-\sqrt{-1}f^0_{2(pn+n)})b_{p,n,0}\eta^{-(p+1)n(j-1)}\\
&\qquad +(f^0_{2(pn+n)-1}+\sqrt{-1}f^0_{2(pn+n)})b_{p,n,0}\eta^{(p+1)n(j-1)}\Big]\\
&= \sqrt{\frac{n}{2}}\sum_{g=0}^{s-1}\Big[\sum_{k=0}^{\infty}(f^0_{2(kN+gn+n)-1}-\sqrt{-1}f^0_{2(kN+gn+n)})b_{g,n,0}\eta^{-(gn+n)(j-1)}\\
&\qquad +\sum_{k=0}^{\infty}(f^0_{2(kN+gn+n)-1}+\sqrt{-1}f^0_{2(kN+gn+n)})b_{g,n,0}\eta^{(gn+n)(j-1)}\Big]\\
&= \sqrt{\frac{n}{2}}\sum_{g=0}^{s-1}\Big[C_{g,0}b_{g,n,0}\eta^{-(gn+n)(j-1)}+D_{g,n}b_{g,n,0}\eta^{(gn+n)(j-1)}\Big], \tag{6.28}
\end{aligned}
$$

which, together with Cauchy-Schwartz inequality, (6.22)–(6.23), and the trivial fact $|b_{g,n,0}|\leq 2$ for $0\leq g\leq s-1$, leads to

$$
\begin{aligned}
\sum_{j=1}^{s}|\delta_{j,0}|^2 &\leq n\sum_{j=1}^{s}\big|\sum_{g=0}^{s-1}C_{g,0}b_{g,n,0}\eta^{-(gn+n)(j-1)}\big|^2+n\sum_{j=1}^{s}\big|\sum_{g=0}^{s-1}D_{g,n}b_{g,n,0}\eta^{(gn+n)(j-1)}\big|^2\\
&= N\left(\sum_{g=0}^{s-1}|C_{g,0}|^2b_{g,n,0}^2+\sum_{g=0}^{s-1}|D_{g,n}|^2b_{g,n,0}^2\right)\\
&\leq 2c_m'n^{-2m}N\sum_{g=0}^{s-1}\sum_{k=0}^{\infty}(|f^0_{2(kN+gn+n)-1}|^2+|f^0_{2(kN+gn+n)}|^2)\times(2\pi(kN+gn+n))^{2m}. \tag{6.29}
\end{aligned}
$$

Combining (6.27) and (6.29) we get that

$$
\begin{aligned}
\sum_{j=1}^{s}\sum_{i=1}^{n}(E\{\bar{f}\}(t_{i,j})-f_0(t_{i,j}))^2 &= \sum_{r=0}^{n-1}\sum_{g=0}^{s-1}|\delta_{j,r}|^2\\
&\leq 2c_m'(\lambda+n^{-2m})N\sum_{i=1}^{n}\sum_{g=0}^{s-1}\sum_{k=0}^{\infty}(|f^0_{2(kN+gn+i)}|^2+|f^0_{2(kN+gn+i)-1}|^2)\\
&\qquad \times(2\pi(kN+gn+i))^{2m}=2c_m'(\lambda+n^{-2m})NJ(f_0). \tag{6.30}
\end{aligned}
$$

Next we will apply (6.30) to show (3.2). Since $\widehat{f}_j$ is the minimizer of $\ell_{j,n,\lambda}(f)$, it satisfies for $1 \le j \le s$,

$$-\frac{1}{n}\sum_{i=1}^{n}(Y_{i,j} - \widehat{f}_j(t_{i,j}))K_{t_{i,j}} + \lambda\widehat{f}_j = 0.$$

Taking expectations, we get that

$$\frac{1}{n}\sum_{i=1}^{n}(E\{\widehat{f}_j\}(t_{i,j}) - f_0(t_{i,j}))K_{t_{i,j}} + \lambda E\{\widehat{f}_j\},$$

therefore, $E\{\widehat{f}_j\}$ is the minimizer to the following functional

$$\ell_{0j}(f) = \frac{1}{2n}\sum_{i=1}^{n}(f(t_{i,j}) - f_0(t_{i,j}))^2 + \frac{\lambda}{2}J(f).$$

Define $g_j = E\{\widehat{f}_j\}$. Since $\ell_{0j}(g_j) \le \ell_{0j}(f_0)$, we get

$$\frac{1}{2n}\sum_{i=1}^{n}(g_j(t_{i,j}) - f_0(t_{i,j}))^2 + \frac{\lambda}{2}J(g_j) \le \frac{\lambda}{2}J(f_0).$$

This means that $J(g_j) \le J(f_0)$, leading to

$$(6.31)\qquad \|\frac{1}{s}\sum_{j=1}^{s}g_j^{(m)}\|_2 \le \frac{1}{s}\sum_{j=1}^{s}\|g_j^{(m)}\|_2 \le \sqrt{J(f_0)}.$$

Note that $E\{\bar{f}\} = \frac{1}{s}\sum_{j=1}^{s}g_j$. Define $g(t) = (E\{\bar{f}\}(t) - f_0(t))^2$. By [4, Lemma (2.24), pp. 58], (6.31) and $m \ge 1$ we get that

$$
\begin{aligned}
\left|\frac{1}{N}\sum_{l=0}^{N-1}g(l/N) - \int_0^1 g(t)dt\right| &\le \frac{2}{N}\int_0^1 \left|\frac{1}{s}\sum_{j=1}^{s}g_j(t) - f_0(t)\right| \times \left|\frac{1}{s}\sum_{j=1}^{s}g_j'(t) - f_0'(t)\right|dt\\
&\le \frac{2}{N}\|\frac{1}{s}\sum_{j=1}^{s}g_j - f_0\|_2 \times \|\frac{1}{s}\sum_{j=1}^{s}g_j' - f_0'\|_2\\
&\le \frac{2}{N}\|\frac{1}{s}\sum_{j=1}^{s}g_j^{(m)} - f_0^{(m)}\|_2^2 \le \frac{8J(f_0)}{N}.
\end{aligned}
$$
(6.32)

Combining (6.30) and (6.32) we get that

$$\|E\{\bar{f}\} - f_0\|_2^2 \le c_m^2 J(f_0)(\lambda + n^{-2m} + N^{-1}),$$

where $c_m^2 = \max\{8, 2c_m'\}$. This completes the proof of (3.2). $\qquad\square$

PROOF OF THEOREM 3.2. Suppose $f_0 = \sum_{\nu=1}^{\infty}f_\nu^0\varphi_\nu$ with $f_\nu^0$ satisfying

$$(6.33)\qquad |f_\nu^0|^2 = \begin{cases} Cn^{-1}(2\pi(n+r))^{-2m}, & \nu = 2(n+r)-1, 1 \le r \le n/2,\\ 0, & \text{otherwise.} \end{cases}$$

It is easy to see that $J(f_0) = \sum_{1 \le r \le n/2} |f_{2(n+r)-1}^0|^2 (2\pi(n+r))^{2m} \le C.$

Consider the decomposition (6.19) and let $\delta_{j,r}$ be defined as in (6.20) and (6.28). It can be easily checked that $C_{g,r} = 0$ for $1 \le r \le n/2$ and $0 \le g \le s-1$. Furthermore, for $1 \le r \le n/2$,

$$
\begin{aligned}
\lambda_{c,r} - a_{1,r} &= \sum_{u=0}^{\infty}(2\pi(un+r))^{-2m} + \sum_{u=1}^{\infty}(2\pi(un-r))^{-2m} - \sum_{u=0}^{\infty}(2\pi(uN+n+r))^{-2m} \\
&\quad - \sum_{u=1}^{\infty}(2\pi(uN-n-r))^{-2m} \ge (2\pi r)^{-2m}.
\end{aligned}
$$

Therefore,

$$
\begin{aligned}
b_{1,r,r}^2 &= \left(\frac{\lambda + \lambda_{c,r} - a_{1,r}}{\lambda + \lambda_{c,r}}\right)^2 \\
&\ge \left(\frac{\lambda + (2\pi r)^{-2m}}{\lambda + 2(1 + \bar{c}_m)(2\pi r)^{-2m}}\right)^2 \ge \frac{1}{4(1 + \bar{c}_m)^2}.
\end{aligned}
\tag{6.34}
$$

Using (6.21) and (6.34), we have

$$
\begin{aligned}
\sum_{j=1}^{s}(\bar{\boldsymbol{f}}_j - \mathbf{f}_{0,j})^T(\bar{\boldsymbol{f}}_j - \mathbf{f}_{0,j}) &= \sum_{j=1}^{s}\sum_{r=0}^{n-1}|\delta_{j,r}|^2 \\
&\ge \sum_{1 \le r \le n/2}\sum_{j=1}^{s}|\delta_{j,r}|^2 \\
&= \sum_{1 \le r \le n/2}\frac{N}{2}\sum_{g=0}^{s-1}|C_{g,r}b_{g,n-r,r} + D_{s-1-g,r}b_{s-1-g,r,r}|^2 \\
&= \sum_{1 \le r \le n/2}\frac{N}{2}\sum_{g=0}^{s-1}|D_{s-1-g,r}|^2 b_{s-1-g,r,r}^2 \\
&= \sum_{1 \le r \le n/2}\frac{N}{2}\sum_{g=0}^{s-1}|D_{g,r}|^2 b_{g,r,r}^2 \\
&\ge \sum_{1 \le r \le n/2}\frac{N}{2}|D_{1,r}|^2 b_{1,r,r}^2 \\
&\ge \frac{N}{8(1 + \bar{c}_m)^2}\sum_{1 \le r \le n/2}|f_{2(n+r)-1}^0|^2 \\
&\ge \frac{NC}{16(3\pi)^{2m}(1 + \bar{c}_m)^2}n^{-2m} \equiv a_m NCn^{-2m},
\end{aligned}
$$

where $a_m = \frac{1}{16(3\pi)^{2m}(1 + \bar{c}_m)^2} < 1$ is an absolute constant depending on $m$ only. Then the conclusion follows by (6.32). Proof is completed. □

6.2. *Proofs in Section 4*.

PROOF OF THEOREM 4.1. For $1 \leq j, l \leq s$, define

$$
\begin{aligned}
\Omega_{j,l} & = \frac{1}{n} \sum_{\nu=1}^{\infty} \frac{\Phi_{\nu,j}^T \Phi_{\nu,l}}{\gamma_\nu^2}, \\
\widetilde{\sigma}_{j,l,r} & = \frac{2}{n} \sum_{k=1}^{\infty} \frac{\cos\left(2\pi k \left(\frac{r}{n} - \frac{j-l}{N}\right)\right)}{(2\pi k)^{4m}}, \ r = 0, 1, \ldots, n-1.
\end{aligned}
$$

Clearly $\Omega_{j,l}$ is a circulant matrix with elements $\widetilde{\sigma}_{j,l,0}, \widetilde{\sigma}_{j,l,1}, \ldots, \widetilde{\sigma}_{j,l,n-1}$. Furthermore, by arguments (6.9)–(6.11) we get that

$$
\Omega_{j,l} = M \Gamma_{j,l} M^*, \tag{6.35}
$$

where $M$ is the same as in (6.4), and $\Gamma_{j,l} = \text{diag}(\delta_{j,l,0}, \delta_{j,l,1}, \ldots, \delta_{j,l,n-1})$, with $\delta_{j,l,r}$, for $r = 1, \ldots, n-1$, given by the following

$$
\delta_{j,l,r} = \sum_{q=1}^{\infty} \frac{\eta^{-(qn-r)(j-l)}}{[2\pi(qn-r)]^{4m}} + \sum_{q=0}^{\infty} \frac{\eta^{(qn+r)(j-l)}}{[2\pi(qn+r)]^{4m}}, \tag{6.36}
$$

and for $r = 0$, given by

$$
\delta_{j,l,0} = \sum_{q=1}^{\infty} \frac{\eta^{qn(j-l)} + \eta^{-qn(j-l)}}{(2\pi qn)^{4m}}. \tag{6.37}
$$

Define $A = \text{diag}(\underbrace{(\Sigma + \lambda I_n)^{-1}, \ldots, (\Sigma + \lambda I_n)^{-1}}_{s})$ and

$$
B = \begin{pmatrix}
\Omega_{1,1} & \Omega_{1,2} & \cdots & \Omega_{1,s} \\
\Omega_{2,1} & \Omega_{2,2} & \cdots & \Omega_{2,s} \\
\cdots & \cdots & \cdots & \cdots \\
\Omega_{s,1} & \Omega_{s,2} & \cdots & \Omega_{s,s}
\end{pmatrix}.
$$

Note that $B$ is $N \times N$ symmetric. Under $H_0$, it can be shown that

$$
\begin{aligned}
\|\bar{f}\|_2^2 & = \sum_{\nu=1}^{\infty} \left(\frac{\sum_{l=1}^{s} \Phi_{\nu,l} (\Sigma + \lambda I_n)^{-1} \boldsymbol{\epsilon}_l}{N \gamma_\nu^2}\right)^2 \\
& = \frac{1}{Ns} \sum_{j,l=1}^{s} \boldsymbol{\epsilon}_j^T (\Sigma + \lambda I_n)^{-1} \left(\frac{1}{n} \sum_{\nu=1}^{\infty} \frac{\Phi_{\nu,j}^T \Phi_{\nu,l}}{\gamma_\nu^2}\right) (\Sigma + \lambda I_n)^{-1} \boldsymbol{\epsilon}_l \\
& = \frac{1}{Ns} \sum_{j,l=1}^{s} \boldsymbol{\epsilon}_j^T (\Sigma + \lambda I_n)^{-1} \Omega_{j,l} (\Sigma + \lambda I_n)^{-1} \boldsymbol{\epsilon}_l \\
& = \frac{1}{Ns} \boldsymbol{\epsilon}^T A B A \boldsymbol{\epsilon} = \frac{1}{Ns} \boldsymbol{\epsilon}^T \Delta \boldsymbol{\epsilon},
\end{aligned}
$$

where $\boldsymbol{\epsilon} = (\boldsymbol{\epsilon}_1^T, \ldots, \boldsymbol{\epsilon}_s^T)^T$ and $\Delta \equiv ABA$.

This implies that $T_{N,\lambda} = \boldsymbol{\epsilon}^T \Delta \boldsymbol{\epsilon}/(Ns)$ with $\mu_{N,\lambda} = \text{trace}(\Delta)/(Ns)$ and $\sigma_{N,\lambda}^2 = 2\text{trace}(\Delta^2)/(Ns)^2$. Define $U = (T_{N,\lambda} - \mu_{N,\lambda})/\sigma_{N,\lambda}$. Then for any $t \in (-1/2, 1/2)$,

$$
\begin{aligned}
\log E\{\exp(tU)\} &= \log E\{\exp(t\boldsymbol{\epsilon}^T \Delta \boldsymbol{\epsilon}/(Ns\sigma_{N,\lambda}))\} - t\mu_{N,\lambda}/\sigma_{N,\lambda} \\
&= -\frac{1}{2}\log\det(I_N - 2t\Delta/(Ns\sigma_{N,\lambda})) - t\mu_{N,\lambda}/\sigma_{N,\lambda} \\
&= t\text{trace}(\Delta)/(Ns\sigma_{N,\lambda}) + t^2\text{trace}(\Delta^2)/((Ns)^2\sigma_{N,\lambda}^2) \\
&\quad + O(t^3\text{trace}(\Delta^3)/((Ns)^3\sigma_{N,\lambda}^3)) - t\mu_{N,\lambda}/\sigma_{N,\lambda} \\
&= t^2/2 + O(t^3\text{trace}(\Delta^3)/((Ns)^3\sigma_{N,\lambda}^3)).
\end{aligned}
$$

It remains to show that $\text{trace}(\Delta^3)/((Ns)^3\sigma_{N,\lambda}^3) = o(1)$ in order to conclude the proof.

In other words, we need to study $\text{trace}(\Delta^2)$ (used in $\sigma_{N,\lambda}^2$) and $\text{trace}(\Delta^3)$. We start from the former. By direct calculations, we get

$$
\begin{aligned}
\text{trace}(\Delta^2) &= \text{trace}(A^2 B A^2 B) \\
&= \sum_{l=1}^{s} \text{trace}\left(\sum_{j=1}^{s} M(\Lambda_c + \lambda I_n)^{-2}\Gamma_{l,j}(\Lambda_c + \lambda I_n)^{-2}\Gamma_{j,l}M^*\right) \\
&= \sum_{j,l=1}^{s} \text{trace}\left((\Lambda_c + \lambda I_n)^{-2}\Gamma_{l,j}(\Lambda_c + \lambda I_n)^{-2}\Gamma_{j,l}\right) = \sum_{j,l=1}^{s}\sum_{r=0}^{n-1} \frac{|\delta_{j,l,r}|^2}{(\lambda + \lambda_{c,r})^4}.
\end{aligned}
$$

For $1 \leq g \leq s$ and $0 \leq r \leq n-1$, define

$$
A_{g,r} = \sum_{p=0}^{\infty} \frac{1}{[2\pi(pN + gn - r)]^{4m}}.
$$

Using (6.36) and (6.37), it can be shown that for $r = 1, 2, \ldots, n-1$,

$$
\begin{aligned}
\sum_{j,l=1}^{s} |\delta_{j,l,r}|^2 &= \sum_{j,l=1}^{s}\left|\sum_{g=1}^{s} A_{g,r}\eta^{-gn(j-l)} + \sum_{g=1}^{s} A_{g,n-r}\eta^{(g-1)n(j-l)}\right|^2 \\
&= \sum_{j,l=1}^{s}\left(\sum_{g,g'=1}^{s} A_{g,r}A_{g',r}\eta^{-(g-g')n(j-l)} + \sum_{g,g'=1}^{s} A_{g,n-r}A_{g',n-r}\eta^{(g-g')n(j-l)}\right. \\
&\quad \left. + \sum_{g,g'=1}^{s} A_{g,r}A_{g',n-r}\eta^{-(g+g'-1)n(j-l)} + \sum_{g,g'=1}^{s} A_{g,n-r}A_{g',r}\eta^{(g+g'-1)n(j-l)}\right) \\
&= s^2\sum_{g=1}^{s} A_{g,r}^2 + s^2\sum_{g=1}^{s} A_{g,n-r}^2 + 2s^2\sum_{g=1}^{s} A_{g,r}A_{s+1-g,n-r} \\
&\geq s^2\sum_{g=1}^{s} A_{g,r}^2 + s^2\sum_{g=1}^{s} A_{g,n-r}^2.
\end{aligned}
$$

(6.38)

Since

$$
(6.39) \qquad \sum_{g=1}^{s} A_{g,r}^2 = \sum_{g=1}^{s}\left(\sum_{p=0}^{\infty}\frac{1}{[2\pi(pN + gn - r)]^{4m}}\right)^2 \geq \frac{1}{[2\pi(n-r)]^{8m}},
$$

we get that

$$
\begin{aligned}
\mathrm{trace}(\Delta^2) \quad &\geq \quad \sum_{r=1}^{n-1} \frac{s^2(\sum_{g=1}^s A_{g,r}^2 + \sum_{g=1}^s A_{g,n-r}^2)}{(\lambda + \lambda_{c,r})^4} \\
&\geq \quad s^2 \sum_{r=1}^{n-1} \frac{\frac{1}{[2\pi(n-r)]^{8m}} + \frac{1}{(2\pi r)^{8m}}}{(\lambda + \lambda_{c,r})^4} \\
&\geq \quad \frac{2s^2}{(2 + 2\bar{c}_m)^4} \sum_{1 \leq r \leq n/2} \frac{\frac{1}{(2\pi r)^{8m}}}{(\lambda + \frac{1}{(2\pi r)^{2m}})^4} \\
&= \quad \frac{s^2}{8(1 + \bar{c}_m)^4} \sum_{1 \leq r \leq n/2} \frac{1}{(1 + (2\pi r h)^{2m})^4} \\
&\geq \quad \frac{s^2}{8(1 + \bar{c}_m)^4} h^{-1} \int_h^{nh/2} \frac{1}{(1 + (2\pi x)^{2m})^4} dx.
\end{aligned}
$$

Meanwhile, (6.38) indicates that for $1 \leq r \leq n - 1$,

$$
\sum_{j,l=1}^s |\delta_{j,l,r}|^2 \leq 2s^2 \sum_{g=1}^s A_{g,r}^2 + 2s^2 \sum_{g=1}^s A_{g,n-r}^2.
$$

From (6.39) we get that for $1 \leq r \leq n - 1$,

$$
\sum_{g=1}^s A_{g,r}^2 \leq \frac{c_m}{(2\pi(n-r))^{8m}},
$$

where $c_m > 0$ is a constant depending on $m$ only.

Similar analysis to (6.38) shows that

$$
\begin{aligned}
\sum_{j,l=1}^s |\delta_{j,l,0}|^2 \quad &= \quad \sum_{j,l=1}^s \left| \sum_{g=1}^s A_{g,0}(\eta^{gn(j-l)} + \eta^{-gn(j-l)}) \right|^2 \\
&= \quad 2s^2 \sum_{g=1}^s A_{g,0}^2 + 2s^2 \sum_{g=1}^{s-1} A_{g,0} A_{s-g,0} + 2s^2 A_{s,0}^2 \\
&\leq \quad 4s^2 \sum_{g=1}^s A_{g,0}^2 \leq c_m s^2 (2\pi n)^{-8m}.
\end{aligned}
$$

Therefore,

$$
\begin{aligned}
\mathrm{trace}(\Delta^2) \quad &\leq \quad \frac{4s^2 \sum_{g=1}^s A_{g,0}^2}{(\lambda + \lambda_{c,0})^4} + 2s^2 \sum_{r=1}^{n-1} \frac{\sum_{g=1}^s A_{g,r}^2 + \sum_{g=1}^s A_{g,n-r}^2}{(\lambda + \lambda_{c,r})^4} \\
&\leq \quad 4c_m s^2 \sum_{r=1}^n \frac{1}{(1 + (2\pi r h)^{2m})^4} \\
&\leq \quad 4c_m s^2 h^{-1} \int_0^{nh} \frac{1}{(1 + (2\pi x)^{2m})^4} dx.
\end{aligned}
$$

By the above statements, we get that

$$(6.40) \qquad \sigma_{N,\lambda}^2 = 2\mathrm{trace}(\Delta^2)/(Ns)^2 \asymp \begin{cases} \frac{n}{N^2}, & \text{if } nh \to 0, \\ \frac{1}{N^2 h}, & \text{if } \lim_N nh > 0. \end{cases}$$

To the end, we look at the trace of $\Delta^3$. By direct examinations, we have

$$\begin{aligned}
\mathrm{trace}(\Delta^3) &= \mathrm{trace}(ABA^2BA^2BA) \\
&= \sum_{j,k=1}^{s} \mathrm{trace}\left( [\sum_{l=1}^{s} M(\Lambda_c + \lambda I_n)^{-2}\Gamma_{j,l}(\Lambda_c + \lambda I_n)^{-2}\Gamma_{l,k}M^*] \right. \\
&\qquad\qquad \left. \times M(\Lambda_c + \lambda I_n)^{-2}\Gamma_{k,j}M^* \right) \\
&= \sum_{j,k,l=1}^{s} \mathrm{trace}\left( (\Lambda_c + \lambda I_n)^{-2}\Gamma_{j,l}(\Lambda_c + \lambda I_n)^{-2}\Gamma_{l,k}(\Lambda_c + \lambda I_n)^{-2}\Gamma_{k,j} \right) \\
&= \sum_{j,k,l=1}^{s} \sum_{r=0}^{n-1} \frac{\delta_{j,l,r}\delta_{l,k,r}\delta_{k,j,r}}{(\lambda + \lambda_{c,r})^6}.
\end{aligned}$$

For $r = 1, 2, \ldots, n-1$, it can be shown that

$$\begin{aligned}
(6.41) \quad \delta_{j,l,r}\delta_{l,k,r}\delta_{k,j,r} &= \left( \sum_{q=1}^{\infty} \frac{\eta^{-qn(j-l)}}{(2\pi(qn-r))^{4m}} + \sum_{q=0}^{\infty} \frac{\eta^{qn(j-l)}}{(2\pi(qn+r))^{4m}} \right) \\
&\quad \times \left( \sum_{q=1}^{\infty} \frac{\eta^{-qn(l-k)}}{(2\pi(qn-r))^{4m}} + \sum_{q=0}^{\infty} \frac{\eta^{qn(l-k)}}{(2\pi(qn+r))^{4m}} \right) \\
&\quad \times \left( \sum_{q=1}^{\infty} \frac{\eta^{-qn(k-j)}}{(2\pi(qn-r))^{4m}} + \sum_{q=0}^{\infty} \frac{\eta^{qn(k-j)}}{(2\pi(qn+r))^{4m}} \right).
\end{aligned}$$

We next proceed to show that for $1 \le r \le n-1$,

$$(6.42) \quad \sum_{l,j,k=1}^{s} \delta_{j,l,r}\delta_{l,k,r}\delta_{k,j,r} \le \frac{96m}{12m-1}\left(\frac{4m}{4m-1}\right)^3 s^3\left(\frac{1}{(2\pi(n-r))^{12m}} + \frac{1}{(2\pi r)^{12m}}\right).$$

Using the trivial fact that $A_{g,r} \leq \frac{4m}{4m-1} \times \frac{1}{(2\pi(gn-r))^{4m}}$, the first term in (6.41) satisfies

$$\sum_{j,l,k=1}^{s} \sum_{q_1=1}^{\infty} \frac{\eta^{-q_1 n(j-l)}}{(2\pi(q_1 n - r))^{4m}} \sum_{q_2=1}^{\infty} \frac{\eta^{-q_2 n(j-l)}}{(2\pi(q_2 n - r))^{4m}} \sum_{q_3=1}^{\infty} \frac{\eta^{-q_3 n(j-l)}}{(2\pi(q_3 n - r))^{4m}}$$

$$= \sum_{j,l,k=1}^{s} \sum_{g_1=1}^{s} A_{g_1,r} \eta^{-g_1 n(j-l)} \sum_{g_2=1}^{s} A_{g_2,r} \eta^{-g_2 n(l-k)} \sum_{g_3=1}^{s} A_{g_3,r} \eta^{-g_3 n(k-j)}$$

$$= \sum_{g_1,g_2,g_3=1}^{s} A_{g_1,r} A_{g_2,r} A_{g_3,r} \sum_{j,l,k=1}^{s} \eta^{-g_1 n(j-l)} \eta^{-g_2 n(l-k)} \eta^{-g_3 n(k-j)}$$

$$= \sum_{g_1,g_2,g_3=1}^{s} A_{g_1,r} A_{g_2,r} A_{g_3,r} \sum_{j=1}^{s} \eta^{(g_3-g_1)n(j-1)} \sum_{l=1}^{s} \eta^{(g_1-g_2)n(l-1)} \sum_{k=1}^{s} \eta^{(g_2-g_3)n(k-1)}$$

$$= s^3 \sum_{g=1}^{s} A_{g,r}^3$$

$$\leq \left(\frac{4m}{4m-1}\right)^3 s^3 \sum_{g=1}^{s} \frac{1}{(2\pi(gn-r))^{12m}}$$

$$\leq \frac{12m}{12m-1} \left(\frac{4m}{4m-1}\right)^3 s^3 \frac{1}{(2\pi(n-r))^{12m}}.$$

Similarly, one can show that all other terms in (6.41) are upper bounded by

$$\frac{12m}{12m-1} \left(\frac{4m}{4m-1}\right)^3 s^3 \left(\frac{1}{(2\pi(n-r))^{12m}} + \frac{1}{(2\pi r)^{12m}}\right).$$

Therefore, (6.42) holds. It can also be shown by (6.37) and similar analysis that

$$(6.43) \qquad \sum_{j,l,k=1}^{s} \delta_{j,l,0} \delta_{l,k,0} \delta_{k,j,0} \leq s^3 (2\pi n)^{-12m}.$$

Using (6.42) and (6.43), one can get that

$$\text{trace}(\Delta^3) = \sum_{l,j,k=1}^{s} \sum_{r=0}^{n-1} \frac{\delta_{j,l,r} \delta_{l,k,r} \delta_{k,j,r}}{(\lambda + \lambda_{c,r})^6}$$

$$\lesssim s^3 \sum_{r=1}^{n-1} \frac{\frac{1}{(2\pi(n-r))^{12m}} + \frac{1}{(2\pi r)^{12m}}}{(\lambda + \lambda_{c,r})} + s^3 \frac{\frac{1}{(2\pi n)^{12m}}}{(\lambda + \lambda_{c,0})^{12m}}$$

$$\lesssim s^3 \sum_{r=1}^{n} \frac{1}{(1 + (2\pi rh)^{2m})^6}$$

$$(6.44) \qquad \lesssim s^3 h^{-1} \int_0^{nh} \frac{1}{(1 + (2\pi x)^{2m})^6} dx \asymp \begin{cases} s^3 n, & \text{if } nh \to 0, \\ s^3 h^{-1}, & \text{if } \lim_N nh > 0. \end{cases}$$

Combining (6.40) and (6.44), and using the assumptions $n \to \infty$, $h \to 0$, we get that

$$\text{trace}(\Delta^3)/((Ns)^3 \sigma_{N,\lambda}^3) \lesssim \begin{cases} n^{-1/2}, & \text{if } nh \to 0, \\ h^{1/2}, & \text{if } \lim_N nh > 0. \end{cases} = o(1).$$

Proof is completed.                                                                                 □

PROOF OF THEOREM 4.2. Throughout the proof, we assume that data $Y_1, \ldots, Y_N$ are generated from the sequence of alternative hypotheses: $f \in \mathcal{B}$ and $\|f\|_2 \geq C_\varepsilon d_{N,\lambda}$. Define $\mathbf{f}_j = (f(t_{1,j}), \ldots, f(t_{n,j}))^T$ for $1 \leq j \leq s$. Then it can be shown that

$$
\begin{aligned}
N s T_{N,\lambda} &= N s \sum_{\nu=1}^{\infty} \bar{f}_\nu^2 \\
&= \sum_{j,l=1}^{s} Y_j^T (\Sigma + \lambda I_n)^{-1} \Omega_{j,l} (\Sigma + \lambda I_n)^{-1} Y_l \\
&= \sum_{j,l=1}^{s} Y_j^T M (\Lambda_c + \lambda I_n)^{-1} \Gamma_{j,l} (\Lambda_c + \lambda I_n)^{-1} M^* Y_l \\
&= \sum_{j,l=1}^{s} \mathbf{f}_j^T M (\Lambda_c + \lambda I_n)^{-1} \Gamma_{j,l} (\Lambda_c + \lambda I_n)^{-1} M^* \mathbf{f}_l \\
&\quad + \sum_{j,l=1}^{s} \mathbf{f}_j^T M (\Lambda_c + \lambda I_n)^{-1} \Gamma_{j,l} (\Lambda_c + \lambda I_n)^{-1} M^* \boldsymbol{\epsilon}_l \\
&\quad + \sum_{j,l=1}^{s} \boldsymbol{\epsilon}_j^T M (\Lambda_c + \lambda I_n)^{-1} \Gamma_{j,l} (\Lambda_c + \lambda I_n)^{-1} M^* \mathbf{f}_l \\
&\quad + \sum_{j,l=1}^{s} \boldsymbol{\epsilon}_j^T M (\Lambda_c + \lambda I_n)^{-1} \Gamma_{j,l} (\Lambda_c + \lambda I_n)^{-1} M^* \boldsymbol{\epsilon}_l \\
&\equiv T_1 + T_2 + T_3 + T_4.
\end{aligned}
\tag{6.45}
$$

Next we will analyze all the four terms in the above. Let $f = \sum_{\nu=1}^{\infty} f_\nu \varphi_\nu$. For $0 \leq r \leq n-1$ and $1 \leq l \leq s$, define $d_{l,r} = x_r^* \mathbf{f}_l$. Then it holds that

$$
d_{l,r} = \sum_{p=0}^{\infty} \sum_{v=1}^{n} f_{2(pn+v)-1} x_r^* \Phi_{2(pn+v)-1,l}^T + \sum_{p=0}^{\infty} \sum_{v=1}^{n} f_{2(pn+v)} x_r^* \Phi_{2(pn+v),l}^T.
$$

Using (6.12) and (6.13), we get that for $1 \leq r \leq n-1$,

$$
\begin{aligned}
d_{l,r} &= \sum_{p=0}^{\infty} \sum_{v=1}^{n-1} f_{2(pn+v)-1} \left( \sqrt{\frac{n}{2}} \left( \eta^{-(pn+v)(l-1)} I(r+v=n) + \eta^{(pn+v)(l-1)} I(r=v) \right) \right. \\
&\quad + \sum_{p=0}^{\infty} \sum_{v=1}^{n-1} f_{2(pn+v)} \left( -\sqrt{\frac{n}{2}} \right) \left( \eta^{-(pn+v)(l-1)} I(r+v=n) - \eta^{(pn+v)(l-1)} I(r=v) \right) \\
&= \sqrt{\frac{n}{2}} \sum_{p=0}^{\infty} \left[ (f_{2(pn+n-r)-1} - \sqrt{-1} f_{2(pn+n-r)}) \eta^{-(pn+n-r)(l-1)} \right. \\
&\quad + (f_{2(pn+r)-1} + \sqrt{-1} f_{2(pn+r)}) \eta^{(pn+r)(l-1)} \Big],
\end{aligned}
\tag{6.46}
$$

and for $r = 0$,

$$
\begin{aligned}
d_{l,0} &= \sum_{p=0}^{\infty} f_{2(pn+n)-1} x_0^* \Phi_{2(pn+n)-1,l}^T + \sum_{p=0}^{\infty} f_{2(pn+n)} x_0^* \Phi_{2(pn+n),l}^T \\
&= \sqrt{\frac{n}{2}} \sum_{p=0}^{\infty} \Big[ (f_{2(pn+n)-1} - \sqrt{-1} f_{2(pn+n)}) \eta^{-(pn+n)(l-1)} \\
&\qquad + (f_{2(pn+n)-1} + \sqrt{-1} f_{2(pn+n)}) \eta^{(pn+n)(l-1)} \Big].
\end{aligned}
\tag{6.47}
$$

We first look at $T_1$. It can be examined directly that

$$
\begin{aligned}
T_1 &= \sum_{j,l=1}^{s} (\overline{d_{j,0}}, \ldots, \overline{d_{j,n-1}}) \mathrm{diag}\left( \frac{\delta_{j,l,0}}{(\lambda + \lambda_{c,0})^2}, \ldots, \frac{\delta_{j,l,n-1}}{(\lambda + \lambda_{c,n-1})^2} \right) \times (d_{l,0}, \ldots, d_{l,n-1})^T \\
&= \sum_{r=0}^{n-1} \frac{\sum_{j,l=1}^{s} \delta_{j,l,r} \overline{d_{j,r}} d_{l,r}}{(\lambda + \lambda_{c,r})^2}.
\end{aligned}
\tag{6.48}
$$

Using similar arguments as (6.14)–(6.18), one can show that for $p \geq 0$, $1 \leq v \leq n$, $0 \leq r \leq n-1$ and $1 \leq j \leq s$,

$$
\begin{aligned}
\frac{1}{s} \sum_{l=1}^{s} \delta_{j,l,r} x_r^* \Phi_{2(pn+v)-1,l}^T &= b_{p,v} x_r^* \Phi_{2(pn+v)-1,j}^T, \\
\frac{1}{s} \sum_{l=1}^{s} \delta_{j,l,r} x_r^* \Phi_{2(pn+v),l}^T &= b_{p,v} x_r^* \Phi_{2(pn+v),j}^T,
\end{aligned}
\tag{6.49}
$$

where

$$
b_{p,v} = \begin{cases} \sum_{u \geq -p/s} \frac{1}{(2\pi(uN+pn+v))^{4m}} + \sum_{u \geq (p+1)/s} \frac{1}{(2\pi(uN-pn-v))^{4m}}, & \text{for } 1 \leq v \leq n-1, \\ \sum_{u \geq -p/s} \frac{1}{(2\pi(uN+pn+n))^{4m}} + \sum_{u \geq (p+2)/s} \frac{1}{(2\pi(uN-pn-n))^{4m}}, & \text{for } v = n. \end{cases}
$$

By (6.49), we have

$$
\begin{aligned}
\sum_{j,l=1}^{s} \delta_{j,l,r} \overline{d_{j,r}} d_{l,r} &= \sum_{j=1}^{s} \overline{d_{j,r}} \sum_{l=1}^{s} \delta_{j,l,r} d_{l,r} \\
&= \sum_{j=1}^{s} \overline{d_{j,r}} \sum_{l=1}^{s} \delta_{j,l,r} \left( \sum_{p=0}^{\infty} \sum_{v=1}^{n} f_{2(pn+v)-1} x_r^* \Phi_{2(pn+v)-1,l} \right. \\
&\qquad \left. + \sum_{p=0}^{\infty} \sum_{v=1}^{n} f_{2(pn+v)} x_r^* \Phi_{2(pn+v),l}^T \right) \\
&= \sum_{j=1}^{s} \overline{d_{j,r}} \left( \sum_{p=0}^{\infty} \sum_{v=1}^{n} f_{2(pn+v)-1} \sum_{l=1}^{s} \delta_{j,l,r} x_r^* \Phi_{2(pn+v)-1,l}^T \right. \\
&\qquad \left. + \sum_{p=0}^{\infty} \sum_{v=1}^{n} f_{2(pn+v)} \sum_{l=1}^{s} \delta_{j,l,r} x_r^* \Phi_{2(pn+v),l}^T \right) \\
&= s \sum_{j=1}^{s} \overline{d_{j,r}} \left( \sum_{p=0}^{\infty} \sum_{v=1}^{n} f_{2(pn+v)-1} b_{p,v} x_r^* \Phi_{2(pn+v)-1,j}^T \right. \\
&\qquad \left. + \sum_{p=0}^{\infty} \sum_{v=1}^{n} f_{2(pn+v)} b_{p,v} x_r^* \Phi_{2(pn+v),j}^T \right).
\end{aligned}
$$

It then follows from (6.46) and (6.47), trivial facts $b_{s-1-g,r} = b_{g,n-r}$ and $C_{g,n-r} = \overline{D_{g,r}}$ (both $C_{g,r}$ and $D_{g,r}$ are defined similarly as those in the proof of Theorem 3.1, but with $f_0$ therein replaced by $f$), and direct calculations that for $1 \le r \le n-1$

$$
\begin{aligned}
\sum_{j,l=1}^{s} \delta_{j,l,r} \overline{d_{j,r}} d_{l,r} &= \frac{sn}{2} \sum_{j=1}^{s} \sum_{p=0}^{\infty} \Big[ (f_{2(pn+n-r)-1} + \sqrt{-1} f_{2(pn+n-r)}) \eta^{(pn+n-r)(j-1)} \\
&\quad + (f_{2(pn+r)-1} - \sqrt{-1} f_{2(pn+r)}) \eta^{-(pn+r)(j-1)} \Big] \\
&\quad \times \sum_{p=0}^{\infty} \Big[ (f_{2(pn+n-r)-1} - \sqrt{-1} f_{2(pn+n-r)}) b_{p,n-r} \eta^{-(pn+n-r)(j-1)} \\
&\quad + (f_{2(pn+r)-1} + \sqrt{-1} f_{2(pn+r)}) b_{p,r} \eta^{(pn+r)(j-1)} \Big] \\
&= \frac{N}{2} \sum_{j=1}^{s} \sum_{g=0}^{s-1} \sum_{k=0}^{\infty} \Big[ (f_{2(kN+gn+n-r)-1} + \sqrt{-1} f_{2(kN+gn+n-r)}) \eta^{(gn+n-r)(j-1)} \\
&\quad + (f_{2(kN+gn+r)-1} - \sqrt{-1} f_{2(kN+gn+r)}) \eta^{-(gn+r)(j-1)} \Big] \\
&\quad \times \sum_{g=0}^{s-1} \sum_{k=0}^{\infty} \Big[ (f_{2(kN+gn+n-r)-1} - \sqrt{-1} f_{2(kN+gn+n-r)}) b_{ks+g,n-r} \eta^{-(gn+n-r)(j-1)} \\
&\quad + (f_{2(kN+gn+r)-1} + \sqrt{-1} f_{2(kN+gn+r)}) b_{ks+g,r} \eta^{(gn+r)(j-1)} \Big] \\
&= \frac{N}{2} \sum_{j=1}^{s} \left[ \sum_{g=0}^{s-1} \overline{C_{g,r}} \eta^{(gn+n)(j-1)} + \sum_{g=0}^{s-1} \overline{D_{g,r}} \eta^{-gn(j-1)} \right] \\
&\quad \times \left[ \sum_{g=0}^{s-1} b_{g,n-r} C_{g,r} \eta^{-(gn+n)(j-1)} + \sum_{g=0}^{s-1} b_{g,r} D_{g,r} \eta^{gn(j-1)} \right] \\
&= \frac{Ns}{2} \left( \sum_{g=0}^{s-1} b_{g,n-r} |C_{g,r}|^2 + \sum_{g=0}^{s-1} b_{s-1-g,r} \overline{C_{g,r}} D_{s-1-g,r} \right. \\
&\quad \left. + \sum_{g=0}^{s-1} b_{s-1-g,n-r} \overline{D_{g,r}} C_{s-1-g,r} + \sum_{g=0}^{s-1} b_{g,r} |D_{g,r}|^2 \right),
\end{aligned}
$$

which leads to

$$
\tag{6.50}
\sum_{r=1}^{n-1} \frac{\sum_{j,l=1}^{s} \delta_{j,l,r} \overline{d_{j,r}} d_{l,r}}{(\lambda + \lambda_{c,r})^2} = \frac{Ns}{2} \sum_{r=1}^{n-1} \frac{\sum_{g=0}^{s-1} b_{g,r} |C_{s-1-g,r} + D_{g,r}|^2}{(\lambda + \lambda_{c,r})^2}.
$$

Since $J(f) \le C$, equivalently, $\sum_{\nu=1}^{\infty} (f_{2\nu-1}^2 + f_{2\nu}^2)(2\pi\nu)^{2m} \le C$, we get that

$$
\tag{6.51}
\sum_{1 \le r \le n/2} (f_{2r-1}^2 + f_{2r}^2) \ge \|f\|_2^2 - C(2\pi n)^{-2m}.
$$

Meanwhile, for $1 \le r \le n/2$, using similar arguments as (6.25) and (6.26) one can show that there exists a constant $c'_m$ relying on $C$ and $m$ s.t.

$$
\begin{aligned}
|C_{s-1,r} + D_{0,r}|^2 &= \left( f_{2r-1} + \sum_{k=0}^{\infty} f_{2(kN+N-r)-1} + \sum_{k=1}^{\infty} f_{2(kN+r)-1} \right)^2 \\
&\quad + \left( f_{2r} + \sum_{k=1}^{\infty} f_{2(kN+r)} - \sum_{k=0}^{\infty} f_{2(kN+N-r)} \right)^2 \\
&\geq \frac{1}{2}(f_{2r-1}^2 + f_{2r-1}^2) - c'_m N^{-2m},
\end{aligned}
\tag{6.52}
$$

and

$$
\begin{aligned}
|C_{s-1,r} + D_{0,r}|^2 (2\pi r)^{2m} &\leq 4 \left[ (\sum_{k=0}^{\infty} f_{2(kN+N-r)-1})^2 + (\sum_{k=0}^{\infty} f_{2(kN+N-r)})^2 \right.\\
&\qquad \left. + (\sum_{k=0}^{\infty} f_{2(kN+r)-1})^2 + (\sum_{k=0}^{\infty} f_{2(kN+r)})^2 \right] (2\pi r)^{2m} \\
&\leq \left( 4\sum_{k=0}^{\infty} f_{2(kN+N-r)-1}^2 (2\pi(kN+N-r))^{2m} \sum_{k=0}^{\infty}(2\pi(kN+N-r))^{-2m} \right. \\
&\qquad + 4\sum_{k=0}^{\infty} f_{2(kN+N-r)}^2 (2\pi(kN+N-r))^{2m} \sum_{k=0}^{\infty}(2\pi(kN+N-r))^{-2m} \\
&\qquad + 4\sum_{k=0}^{\infty} f_{2(kN+r)-1}^2 (2\pi(kN+r))^{2m} \sum_{k=0}^{\infty}(2\pi(kN+r))^{-2m} \\
&\qquad \left. + 4\sum_{k=0}^{\infty} f_{2(kN+r)}^2 (2\pi(kN+r))^{2m} \sum_{k=0}^{\infty}(2\pi(kN+r))^{-2m} \right) \times (2\pi r)^{2m} \\
&\leq \left( \frac{8m}{2m-1} \sum_{k=0}^{\infty}(f_{2(kN+N-r)-1}^2 + f_{2(kN+N-r)}^2)\gamma_{kN+N-r}(2\pi(N-r))^{-2m} \right. \\
&\qquad \left. + \frac{8m}{2m-1} \sum_{k=0}^{\infty}(f_{2(kN+r)-1}^2 + f_{2(kN+r)}^2)\gamma_{kN+r}(2\pi r)^{-2m} \right) \times (2\pi r)^{2m},
\end{aligned}
$$

which, together with the fact $N \geq 2r$ for $1 \leq r \leq n/2$, leads to that

$$
\sum_{1 \leq r \leq n/2} |C_{s-1,r} + D_{0,r}|^2 (2\pi r)^{2m} \leq c'_m.
\tag{6.53}
$$

Furthermore, it can be verified that for $1 \leq r \leq n/2$,

$$
\begin{aligned}
\frac{\lambda_{c,r}^2 - b_{0,r}}{(\lambda + \lambda_{c,r})}(2\pi r)^{-2m} &\leq \frac{((2\pi r)^{-2m} + (2\pi(n-r))^{-2m} + \bar{c}_m(2\pi n)^{-2m})^2 - (2\pi r)^{-4m}}{((2\pi r)^{-2m} + (2\pi(n-r))^{-2m})^2}(2\pi r)^{-2m} \\
&\leq c'_m n^{-2m},
\end{aligned}
\tag{6.54}
$$

which leads to that

$$
\begin{aligned}
\frac{(\lambda + \lambda_{c,r})^2 - b_{0,r}}{(\lambda + \lambda_{c,r})^2}(2\pi r)^{-2m} &= \frac{\lambda^2 + 2\lambda\lambda_{c,r}}{(\lambda + \lambda_{c,r})^2}(2\pi r)^{-2m} + \frac{\lambda_{c,r}^2 - b_{0,r}}{(\lambda + \lambda_{c,r})^2}(2\pi r)^{-2m} \\
&\leq 2\lambda + c'_m n^{-2m}.
\end{aligned}
\tag{6.55}
$$

Then, using (6.48)–(6.50) and (6.51)–(6.55) one gets that

$$
\begin{aligned}
T_1 &\geq \frac{Ns}{2} \sum_{1 \leq r \leq n/2} \frac{b_{0,r}|C_{s-1,r} + D_{0,r}|^2}{(\lambda + \lambda_{c,r})^2} \\
&= \frac{Ns}{2} \left[ \sum_{1 \leq r \leq n/2} |C_{s-1,r} + D_{0,r}|^2 - \sum_{1 \leq r \leq n/2} \frac{(\lambda + \lambda_{c,r})^2 - b_{0,r}}{(\lambda + \lambda_{c,r})^2} |C_{s-1,r} + D_{0,r}|^2 \right] \\
&\geq \frac{Ns}{2} \left( \frac{1}{2} \|f\|_2^2 - c_m' n^{-2m} - c_m' N^{-2m} - c_m'(2\lambda + c_m' n^{-2m}) \right) \geq C'Ns\sigma_{N,\lambda},
\end{aligned}
$$

(6.56)

where the last inequality follows by $\|f\|_2^2 \geq 4C'(\lambda + n^{-2m} + \sigma_{N,\lambda})$ for a large constant $C'$ satisfying $2C' > 2c_m' + (c_m')^2$. To achieve the desired power, we need to enlarge $C'$ further. This will be described later. Combining (6.56) with (6.40) and (6.56) we get that

(6.57)                $T_1 \gg s$   uniformly for $f \in \mathcal{B}$ with $\|f\|_2^2 \geq 4C'd_{N,\lambda}^2$.

Terms $T_2$ and $T_3$ can be handled similarly. To handle $T_2$, note that $T_2 = \mathbf{f}^T \Delta \boldsymbol{\epsilon}$, where $\mathbf{f} = (\mathbf{f}_1^T, \ldots, \mathbf{f}_s^T)^T$, $\boldsymbol{\epsilon} = (\boldsymbol{\epsilon}_1^T, \ldots, \boldsymbol{\epsilon}_s^T)^T$, and $\Delta$ is defined in the proof of Theorem 4.1. We need to establish $\Delta \leq sI_N$. Define an arbitrary $\mathbf{a} = (a_1^T, \ldots, a_s^T)^T \in \mathbb{R}^N$, where each $a_j$ is an (real) $n$-vector. Let $\xi_j = M^* a_j$ and $\xi = (\xi_1^T, \ldots, \xi_s^T)^T$. For simplicity, put $\xi_j = (\xi_{j,0}, \ldots, \xi_{j,n-1})^T$ for $1 \leq j \leq s$. Then based on (6.37) and (6.36), we have

$$
\begin{aligned}
\mathbf{a}^T \Delta \mathbf{a} &= \xi^*[(\Lambda_c + \lambda I_n)^{-1} \Gamma_{j,l} (\Lambda_c + \lambda I_n)^{-1}]_{1 \leq j,l \leq s} \xi \\
&= \sum_{r=0}^{n-1} \sum_{j,l=1}^{s} \overline{\xi_{j,r}} \xi_{l,r} \frac{\delta_{j,l,r}}{(\lambda + \lambda_{c,r})^2} \\
&\leq \sum_{r=1}^{n-1} \frac{s \left( \frac{1}{(2\pi(n-r))^{4m}} + \frac{1}{(2\pi r)^{4m}} \right) \sum_{j=1}^{s} |\xi_{j,r}|^2}{(\lambda + \lambda_{c,r})^2} \\
&\quad + \frac{2s \sum_{q=1}^{\infty} \frac{1}{(2\pi qn)^{4m}} \sum_{j=1}^{s} |\xi_{j,0}|^2}{(\lambda + \lambda_{c,0})^2} \\
&\leq s \sum_{r=0}^{n-1} \sum_{j=1}^{s} |\xi_{j,r}|^2 = s\xi^*\xi = s\mathbf{a}^T\mathbf{a},
\end{aligned}
$$

therefore, $\Delta \leq sI_N$. This leads to that, uniformly for $f \in \mathcal{B}$ with $\|f\|_2^2 \geq 4C'd_{N,\lambda}^2$, $E_f\{T_2^2\} = \mathbf{f}^T \Delta^2 \mathbf{f} \leq sT_1$. Together with (6.57), we get that

(6.58)                $\displaystyle\sup_{\substack{f \in \mathcal{B} \\ \|f\|_2^2 \geq 2\sqrt{C'}d_{N,\lambda}}} P_f\left( |T_2| \geq \varepsilon^{-1/2}T_1^{1/2}s^{1/2} \right) \leq \varepsilon.$

Note that (6.58) also applies to $T_3$. By Theorem 4.1, $(T_4/(Ns) - \mu_{N,\lambda})/\sigma_{N,\lambda}$ is $O_P(1)$ uniformly for $f$. Therefore, we can choose $C_\varepsilon' > 0$ s.t. $P_f(|T_4/(Ns) - \mu_{N,\lambda}|/\sigma_{N,\lambda} \geq C_\varepsilon') \leq \varepsilon$ as $N \to \infty$.

It then follows by (6.56), (6.57) and (6.58) that for suitable large $C'$ (e.g., $C' \geq 2(C'_\varepsilon + z_{1-\alpha/2})$), uniformly for $f \in \mathcal{B}$ with $\|f\|_2 \geq 2\sqrt{C'}d_{N,\lambda}$,

$$P_f\left(|T_{N,\lambda} - \mu_{N,\lambda}|/\sigma_{N,\lambda} \geq z_{1-\alpha/2}\right) \leq 3\varepsilon, \text{ as } N \to \infty.$$

Proof is completed. $\qquad\square$

PROOF OF THEOREM 4.3. Define $B_N = \lfloor N^{2/(4m+1)} \rfloor$, the integer part of $N^{2/(4m+1)}$. We prove the theorem in two cases: $\lim_N nh > 0$ and $nh = o(1)$.

**Case I:** $\lim_N nh > 0$.

In this case, it can be shown by $s \gg N^{(4m-1)/(4m+1)}$ (equivalently $n \ll B_N$, leading to $B_N h \gg nh$ hence $B_N h \to \infty$) that $n^{-6m}h^{-4m+1/2}N \ll (B_N/n)^{6m}$. Choose $g$ to be an integer satisfying

$$(6.59) \qquad n^{-6m}h^{-4m+1/2}N \ll g^{6m} \ll (B_N/n)^{6m}.$$

Construct an $f = \sum_{\nu=1}^{\infty} f_\nu \varphi_\nu$ with

$$(6.60) \qquad f_\nu^2 = \begin{cases} \frac{C}{n-1}(2\pi(gn+r))^{-2m}, & \nu = 2(gn+r)-1, \ r = 1, 2, \ldots, n-1, \\ 0, & \text{otherwise.} \end{cases}$$

It can be seen that

$$(6.61) \qquad J(f) = \sum_{r=1}^{s-1} f_{2(gn+r)-1}^2 (2\pi(gn+r))^{2m} = C,$$

and

$$\begin{aligned}
(6.62) \qquad \|f\|_2^2 &= \sum_{r=1}^{n-1} f_{2(gn+r)-1}^2 \\
&= \frac{C}{n-1}\sum_{r=1}^{n-1}(2\pi(gn+r))^{-2m} \\
&\geq C(2\pi(gn+n))^{-2m} = \beta_{N,\lambda}^2 N^{-4m/(4m+1)},
\end{aligned}$$

where $\beta_{N,\lambda}^2 = C[B_N/(2\pi(gn+n))]^{2m}$. Due to (6.59) and $n \ll B_N$, we have $gn+n \ll 2B_N$, which further implies $\beta_{N,\lambda} \to \infty$ as $N \to \infty$.

Using the trivial fact $b_{s-2-g,n} = b_{g,n}$ for $0 \leq g \leq s-2$, one can show that

$$\begin{aligned}
\sum_{j,l=1}^{s} \delta_{j,l,0}\overline{d_{j,0}}d_{l,0} &= \frac{Ns}{2}\left(2\sum_{g'=0}^{s-1}|C_{g',0}|^2 b_{g',n} + \sum_{g'=0}^{s-2}C_{s-2-g',0}C_{g',0}b_{g',n} + C_{s-1,0}^2 b_{s-1,n}\right. \\
&\qquad \left. + \sum_{g'=0}^{s-2}D_{s-2-g',n}D_{g',n}b_{g',n} + D_{s-1,n}^2 b_{s-1,n}\right) \\
(6.63) \qquad &\leq 2Ns\sum_{g'=0}^{s-1}|C_{g',0}|^2 b_{g',n} = 0,
\end{aligned}$$

where the last equality follows by a trivial observation $C_{g',0} = 0$. It follows by (6.63), (6.48) and (6.50) that

$$
\begin{aligned}
T_1 &= \frac{Ns}{2} \sum_{r=1}^{n-1} \frac{\sum_{g'=0}^{s-1} b_{g',r} |C_{s-1-g',r} + D_{g',r}|^2}{(\lambda + \lambda_{c,r})^2} \\
&\leq Ns \sum_{r=1}^{n-1} \frac{\sum_{g'=0}^{s-1} b_{g',r} |C_{s-1-g',r}|^2}{(\lambda + \lambda_{c,r})^2} + Ns \sum_{r=1}^{n-1} \frac{\sum_{g'=0}^{s-1} b_{g',r} |D_{g',r}|^2}{(\lambda + \lambda_{c,r})^2} \\
&= Ns \sum_{r=1}^{n-1} \frac{\sum_{g'=0}^{s-1} b_{s-1-g',n-r} |C_{g',n-r}|^2}{(\lambda + \lambda_{c,r})^2} + Ns \sum_{r=1}^{n-1} \frac{\sum_{g'=0}^{s-1} b_{g',r} |D_{g',r}|^2}{(\lambda + \lambda_{c,r})^2} \\
&= 2Ns \sum_{r=1}^{n-1} \frac{\sum_{g'=0}^{s-1} b_{g',r} |D_{g',r}|^2}{(\lambda + \lambda_{c,r})^2} = 2Ns \sum_{r=1}^{n-1} \frac{b_{g,r} f_{2(gn+r)-1}^2}{(\lambda + \lambda_{c,r})^2},
\end{aligned}
$$

(6.64)

where the last equality follows from the design of $f$, i.e., (6.60). Now it follows from (6.64) and the fact $b_{g,r} \leq c'_m (2\pi(gn + r))^{-4m}$, for some constant $c'_m$ depending on $m$ only, that

$$
\begin{aligned}
T_1 &\leq 2Ns \sum_{r=0}^{n-1} \frac{c'_m (2\pi(gn+r))^{-4m} \frac{C}{n-1} (2\pi(gn+r))^{-2m}}{(\lambda + \lambda_{c,r})^2} \\
&= \frac{2Nsc'_m C}{n-1} \sum_{r=1}^{n-1} \frac{(2\pi(gn+r))^{-6m}}{(\lambda + \lambda_{c,r})^2} \\
&\leq 2c'_m C (2\pi)^{-6m} Ns (gn)^{-6m} h^{-4m} \ll sh^{-1/2} \asymp Ns\sigma_{N,\lambda},
\end{aligned}
$$

(6.65)

where the last "$\asymp$" follows from (4.3).

By (6.58) we have that

$$
|T_2 + T_3| = T_1^{1/2} s^{1/2} O_{P_f}(1) = o_{P_f}(sh^{-1/4}) = o_{P_f}(Ns\sigma_{N,\lambda}).
$$

Hence, by (6.45) and Theorem 4.1 we have

$$
\begin{aligned}
\frac{T_{N,\lambda} - \mu_{N,\lambda}}{\sigma_{N,\lambda}} &= \frac{T_1 + T_2 + T_3}{Ns\sigma_{N,\lambda}} + \frac{T_4/(Ns) - \mu_{N,\lambda}}{\sigma_{N,\lambda}} \\
&= \frac{T_4/(Ns) - \mu_{N,\lambda}}{\sigma_{N,\lambda}} + o_{P_f}(1) \xrightarrow{d} N(0,1).
\end{aligned}
$$

Consequently, as $N \to \infty$

$$
\inf_{\substack{f^\star \in \mathcal{B} \\ \|f^\star\|_2 \geq \beta_{N,\lambda} N^{-2m/(4m+1)}}} P_{f^\star}(\phi_{N,\lambda} = 1) \leq P_f(\phi_{N,\lambda} = 1) \to \alpha.
$$

This shows the desired result in Case I.

**Case II:** $nh = o(1)$.

The proof is similar to Case I although a bit technical difference needs to be emphasized. Since $n \ll B_N$, it can be shown that $Nn^{-2m-1/2} \ll (B_N/n)^{6m}$. Choose $g$ to be an integer satisfying

(6.66)                    $$ Nn^{-2m-1/2} \ll g^{6m} \ll (B_N/n)^{6m}. $$

Let $f = \sum_{\nu=1}^{\infty} f_\nu \varphi_\nu$ with $f_\nu$ satisfying (6.60). Similar to (6.61) and (6.62) one can show that $J(f) = C$ and $\|f\|_2^2 \geq \beta_{N,\lambda}^2 N^{-4m/(4m+1)}$, where $\beta_{N,\lambda}^2 = C[B_N/(2\pi(gn+n))]^{2m}$. It is clear that $\beta_{N,\lambda} \to \infty$ as $N \to \infty$. Then similar to (6.63), (6.48), (6.50) and (6.65) one can show that

$$
\begin{aligned}
T_1 &\leq 2Ns \sum_{r=1}^{n-1} \frac{b_{g,r} f_{2(gn+r)-1}^2}{(\lambda + \lambda_{c,r})^2} \\
&\leq \frac{2Ns c_m' C}{n-1} \sum_{r=1}^{n-1} \frac{(2\pi(gn+r))^{-6m}}{(\lambda + \lambda_{c,r})^2} \\
&\leq 2c_m' C(2\pi)^{-2m} Ns g^{-6m} n^{-2m} \\
&\ll sn^{1/2} \asymp Ns\sigma_{N,\lambda},
\end{aligned}
$$

where the last line follows by (6.66) and (4.3). Then the desired result follows by arguments in the rest of Case I. Proof is completed. □

Department of Mathematical Sciences
Indiana University-Purdue University at Indianapolis
420 University Blvd
Indianapolis, IN 46202
Email: zuofengshang@gmail.com

Department of Statistics
Purdue University
250 N. University Street
West Lafayette, IN 47907
Email: chengg@purdue.edu